\documentclass{template}

\usepackage{tikz}
\usetikzlibrary{arrows.meta,positioning,positioning,patterns,decorations.pathreplacing,shapes.geometric}
\usepackage{marvosym}
\usepackage[eulergreek]{sansmath}
\usepackage{pgfplots}
\usepackage{subcaption}
\usepackage{mathtools}
\usepackage[linesnumbered,ruled,vlined]{algorithm2e}
\usepackage[utf8]{inputenc}
\usepackage[T1]{fontenc}
\usepackage{setspace}
\usepackage{booktabs}
\usepackage[normalem]{ulem}
\usepackage{adjustbox}
\usepackage{enumitem}
\usepackage{amsthm}

\pgfplotsset{compat=1.18}
\ExplSyntaxOn
\ExplSyntaxOff

\definecolor{myLightGray}{RGB}{191,191,191}
\definecolor{myGray}{RGB}{160,160,160}

\tikzstyle{startstop} = [rectangle, rounded corners, 
minimum width=1.5cm, 
minimum height=1cm,
text centered, 
draw=black,
text width=1.5cm]

\tikzstyle{process} = [rectangle, 
minimum width=4cm, 
minimum height=1cm, 
text centered, 
draw=black,
text width=4cm]

\tikzstyle{smallprocess} = [rectangle, 
minimum width=3cm, 
minimum height=1cm, 
text centered, 
draw=black,
text width=3cm]

\tikzstyle{decision} = [diamond, draw, text centered, minimum height=2em, aspect=2, text width=3cm]

\tikzstyle{arrow} = [thick,->,>=stealth]

\allowdisplaybreaks

\begin{document}

\title{Chance-constrained battery management strategies for the electric bus scheduling problem }

\author[1]{L\'ea Ricard$^{*,}$}
\author[2]{Guy Desaulniers}
\author[3]{Andrea Lodi}

\author[2]{Louis-Martin Rousseau}

\affil[1]{School of Architecture, Civil and Environmental Engineering (ENAC), École Polytechnique Fédérale de Lausanne (EPFL), Lausanne, Switzerland}
\affil[2]{Department of Mathematics and Industrial Engineering, Polytechnique Montr\'eal, Montr\'eal, Canada}

\affil[3]{Jacobs Technion-Cornell Institute, Cornell Tech and Technion - IIT, New York, USA}

\date{}

\maketitle

\begingroup
\renewcommand{\thefootnote}{\fnsymbol{footnote}}
\footnotetext[1]{Corresponding author: lea.ricard@epfl.ch}
\endgroup

\begin{abstract}
The global transition to battery electric buses (EBs) presents an opportunity to reduce air and noise pollution in urban areas. However, the adoption of EBs introduces challenges related to limited driving range, extended charging times, and battery degradation. This study addresses these challenges by proposing a novel chance-constrained model for the electric vehicle scheduling problem (E-VSP) that accounts for stochastic energy consumption and battery degradation. The model ensures compliance with recommended state-of-charge (SoC) ranges while optimizing operational costs. A tailored branch-and-price heuristic with stochastic pricing problems is developed. Computational experiments on realistic instances demonstrate that the stochastic approach can provide win-win solutions compared to deterministic baselines in terms of operational costs and battery wear. By limiting the probability of operating EBs outside the recommended SoC range, the proposed framework supports fleet management practices that align with battery leasing company and manufacturer guidelines for battery health and longevity. 
\end{abstract}


\printkeywords

\section{Introduction}

The transition of public transport bus fleets towards battery electric buses (EBs) is gaining momentum worldwide as a solution to reduce local air and noise pollution. Cities like Paris and Copenhagen have already committed to converting their entire bus fleets to EBs by 2025 \citep{Perumal2022}, and it is projected that nearly half of the buses globally will be electric in the near future \citep{Abdelaty2021}. However, the adoption of EBs brings its share of challenges and additional constraints. Compared to traditional internal combustion engine buses, EBs require longer refueling times, have a shorter driving range, and rely on charging infrastructure that is often limited in capacity. Additionally, the batteries used in EBs, predominantly lithium-ion batteries \citep{Zhang2021}, are expensive and represent a significant portion of the cost of ownership. In some cases, the daily cost of battery degradation, for periodical replacement of the batteries, is even more than twice the daily charging cost \citep{Zhou2022}. This replacement is needed due to the progressive loss of capacity during charging and discharging cycles (cyclic aging) and over time during storage (calendar aging). Typically, a battery is considered to have reached its end of life when its capacity declines to 70-80\% of its initial capacity \citep{Lam2012,Zhang2019}.

The degradation mechanisms of lithium-ion batteries are highly complex and difficult to model as they are influenced by various factors, including storage and operating conditions and cell chemistry \citep{Pelletier2017}. Nevertheless, some factors have been identified as accelerators of battery aging, such as overcharging, overdischarging, extreme temperatures, high state-of-charge (SoC) during storage, large depth-of-discharge, high charging and discharging rates, high average SoC, and high SoC deviation \citep{Pelletier2017,Millner2010,Lam2012}. For example, high average SoC and SoC deviation when running EBs in the 0-100\% SoC range (i.e., when batteries are fully recharged and discharged each cycle)  can reduce lithium-ion batteries' expected lifetime by more than 10 years \citep{Lam2012}.\footnote{Expected lifetime estimated using the capacity fading model (measured in kWh/year) for LiFePO$_4$ batteries from \cite{Lam2012}, assuming an EB with a 300 kWh battery capacity, an average energy consumption of 400 kWh per day, and an end of life threshold of 80\%.} To optimize battery life, a safe SoC interval can be applied. A common recommendation is to operate the batteries within the 20-80\% or 30-80\% SoC range most of the time \citep{Kostopoulos2020, Jiang2014, Lu2013}. These safe SoC intervals should be taken into account when planning EB schedules.

At the planning stage, vehicle schedules are determined by solving a vehicle scheduling problem (VSP). This problem consists of assigning vehicles to timetabled trips such that each trip is covered exactly once. A trip is defined by a start time and location, an itinerary composed of a sequence of stops, and an end time and location. A common extension is the multi-depot vehicle scheduling problem (MDVSP), where vehicles can be stored overnight to several depots, and vehicle capacity at each depot is restricted. The VSP and the MDVSP have been extensively studied in the past decades. Several exact and heuristic methods have been proposed, notably those of \cite{Ribeiro1994}, \cite{Lobel1998}, \cite{Hadjar2006}, \cite{Freling2001}, and \cite{Kliewer2006}. A detailed review on the subject is provided in \cite{Desaulniers2007} and \cite{Bunte2010}. 

When managing a fleet of EBs, the VSP must be extended to account for charging activities. This extension, referred to as the electric vehicle scheduling problem (E-VSP), incorporates constraints related to driving range and charger capacity. The E-VSP is proven to be NP-hard, even in the single-depot case \citep{Sassi2017}, and has been widely explored in the literature over the past decade.

\cite{Li2013} is among the first to present a model for the single-depot E-VSP (SD-E-VSP). Their model is suitable for fast-charging or battery-swapping EB fleets, since it assumes that battery service time is fixed. Each charging station is time-expanded, and a time discretization technique is used, so that the capacity of each charging station can be constrained. A column generation-based algorithm is devised, and computational experiments on real-world and randomly generated instances are carried out. In \cite{vanKooten2017}, two alternative SD-E-VSP models are considered. Both models are with continuous and discrete SoC values, respectively. Additionally, the second model incorporates features such as time-of-day electricity pricing, flexible charging processes, and battery wear costs. A column generation-based algorithm is compared with a faster version with Lagrangian relaxation over four datasets. \cite{Wen2016} and \cite{WANG2021} developed an adaptive large neighborhood search heuristic and a genetic column generation-based algorithm, respectively, for solving the multi-depot electric vehicle scheduling problem (MD-E-VSP). \cite{Wu2022} introduced a bi-objective model for the MD-E-VSP considering grid characteristics, namely time-of-use pricing, and peak load risk. The bi-objective problem is reformulated using a lexicographic method and a branch-and-price heuristic featuring a trip chain pool strategy is devised. Furthermore, the E-VSP with a mixed fleet is studied, among others, in \cite{Olsen2020}, \cite{Rinaldi2020}, and \cite{ALVO2021}. The E-VSP integrated with crew scheduling and timetabling is studied in \cite{PERUMAL2021} and \cite{Xu2023}, respectively.

Battery management strategies aimed at optimizing battery life and minimizing wear-and-tear costs have received relatively little attention in the existing literature. To the best of our knowledge, only \cite{Zhang2021} and \cite{Zhou2022} have explicitly addressed the E-VSP while accounting for battery degradation. \cite{Zhang2021} proposed a set partitioning model for the single-terminal E-VSP with nonlinear charging profiles. Each trip is a loop starting and ending at the unique terminal. The total cost of a vehicle schedule is defined as the sum of the charging fees, the cost incurred by battery degradation, and a fixed cost for vehicle acquisition. The capacity fading model of \cite{Lam2012} is employed to estimate the battery degradation cost of a schedule. An exact branch-and-price algorithm with a dual stabilization technique is developed. Their model does not accommodate partial charging, a limitation addressed by \cite{Zhou2022} in their work on the EB charging scheduling problem (EB-CSP). A mixed-integer nonlinear nonconvex program and its mixed-integer linear programming (MILP) approximation are introduced, where the charging and battery degradation functions are represented using linear approximations. Furthermore, \cite{Zheng2022} proposed a model for the EB-CSP with battery degradation under predetermined EB-to-trip assignments. Their model also accounts for the peak-to-average power ratio and time-dependent electricity pricing.

The current literature on the E-VSP with battery degradation suffers from two major shortcomings. First, existing models often rely on capacity fading functions tailored to specific battery types, calibrated using test data that fail to accurately represent real-world driving conditions. This dependence can lead to incorrect conclusions and limits the adaptability of E-VSP models to other battery degradation functions. Second, current research exclusively employs deterministic models, which tend to underestimate the likelihood of critical scenarios such as low SoC, high SoC, or significant SoC deviations---conditions known to accelerate battery aging.

This work aims to fill this gap in the literature by proposing a novel model for the E-VSP that accounts for battery degradation and stochastic energy consumption without relying on a specific capacity fading function. Our model enforces widely used SoC safety ranges, such as 20–80\% or 30–80\%, by limiting the probability of operating an EB outside these ranges. This approach effectively mitigates critical battery degradation factors, including large depth-of-discharge, high SoC, and high SoC deviations, without direct dependence on a specific capacity fading function. To ensure a realistic representation of the system, we incorporate essential E-VSP characteristics, including partial en-route charging, nonlinear charging profiles, and limited-capacity charging stations. Furthermore, we develop a tailored branch-and-price heuristic that features stochastic pricing problems and an exact stochastic dominance rule. 

The remainder of this paper is organized as follows. Section \ref{sec:math_problem} presents a general definition of the E-VSP and introduces our chance-constrained model for the E-VSP with battery degradation and stochastic energy consumption. Furthermore, a method to decompose the chance constraint by vehicle schedule is provided. Section \ref{sec:algo} presents our branch-and-price heuristic. Then, in Section \ref{sec:results}, we compare our approach to deterministic baselines and analyze the tradeoff between the probability of operating within the safe SoC range and operational costs for both approaches. Section \ref{sec:conclusions} summarizes our results.

\section{The E-VSP with battery degradation and stochastic energy consumption}

\label{sec:math_problem}

We first define the E-VSP for the general case, where battery degradation is not considered, in Section \ref{sec:general_def}. In Section \ref{sec:charging_profile}, the implemented non-linear charging function approximation is presented. In Section \ref{sec:CC_model}, we introduce the chance-constrained model for the E-VSP with battery degradation and stochastic energy consumption. A method to compute the probability of complying with the recommended SoC range is provided in Section \ref{sec:proba_overuse}.

\subsection{General case}
\label{sec:general_def}

Let $\mathcal{V}$, $\mathcal{D}$, and $\mathcal{H}$ be a timetable of trips, a set of depots, and a set of charging stations, respectively. For each trip in $\mathcal{V}$, a departure time and location and an arrival time and location are given. Furthermore, the travel time and the energy consumption between any two locations in the EB network---each location being either a starting or ending location of a trip in $\mathcal{V}$, a depot in $\mathcal{D}$, or a charging station in $\mathcal{H}$---is known. The travel time of a trip $i \in V$ is denoted $t_i$. Given this information, the E-VSP consists in finding a set of feasible vehicle schedules $\mathcal{S}^*$ that covers exactly once each timetabled trip while respecting the number of available EBs $b_d$ at each depot $d \in \mathcal{D}$ and the number of available chargers $g_h$ at each charging station $h \in \mathcal{H}$. 

We model the capacity of each charging station with three additional sets. Let  $\mathcal{R}$ be a set of $k + 1$ time intervals, each of duration $\rho$, which partition the planning horizon. Each time interval $r \in R$ is defined by its beginning and end times, denoted by $B_r$ and $E_r$, respectively. Let also $\mathcal{H}^C$ and $\mathcal{H}^W$ be two sets of nodes, where each node is associated with a charging station in $\mathcal{H}$ and a time interval in $\mathcal{R}$. The nodes in $\mathcal{H}^C$, referred to as `charging nodes', represent EBs occupying a charger, while the nodes in $\mathcal{H}^W$, referred to as `waiting nodes', represent EBs waiting at a charging station. For example, we denote by $h^{c1}_{r5} \in \mathcal{H}^C$ and $h^{w1}_{r5}  \in \mathcal{H}^W$ the nodes associated with charging station $h_1 \in \mathcal{H}$ and time interval $r_5 \in \mathcal{R}$, in the charging and the waiting states, respectively.

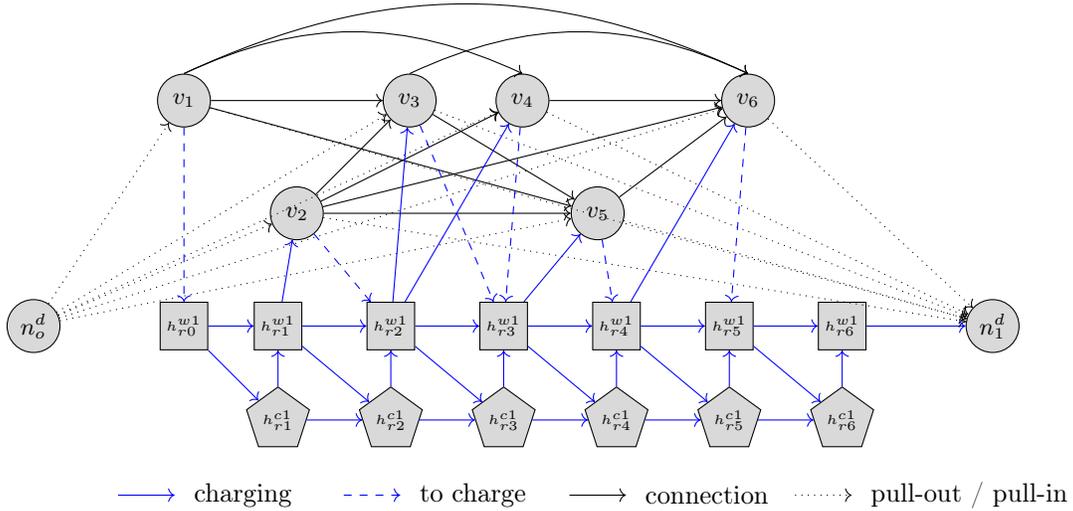
\begin{figure}[!ht]
\centering
\begin{tikzpicture}[
      scale=1, 
      mycircle/.style={
         circle,
         draw=black,
         fill=gray!30,
         fill opacity = 1,
         text opacity=1,
         inner sep=0pt,
         minimum size=20pt,
         font=\small},
      myrectangle/.style={
         rectangle,
         draw=black,
         fill=gray!30,
         fill opacity = 1,
         text opacity=1,
         inner sep=0pt,
         minimum size=18pt,
         font=\tiny},
      mytriangle/.style={
         regular polygon,regular polygon sides=5,
         draw=black,
         fill=gray!30,
         fill opacity = 1,
         text opacity=1,
         inner sep=0pt,
         minimum size=25pt,
         font=\tiny},
      myarrow1/.style={dotted,->, black},
      myarrow2/.style={dashed,->, blue},
      myarrow3/.style={-,->, black},
      myarrow4/.style={-,->, blue},
      node distance=0.6cm and 1.2cm
      ]
      \node[mycircle] at (0,0) (c1) {$n_o^d$};
      \node[mycircle] at (2,3) (c2) {$v_1$};
      \node[mycircle] at (3.5,1.5) (c3) {$v_2$};
      \node[mycircle] at (5,3) (c4) {$v_3$};
      \node[mycircle] at (6.5,3) (c5) {$v_4$};
      \node[mycircle] at (7.5,1.5) (c6) {$v_5$};
      \node[mycircle] at (9.5,3) (c7) {$v_6$};
      \node[mycircle] at (12.75,0) (c8) {$n_1^d$};
      \node[myrectangle] at (2,0) (c9) {$h^{w1}_{r0}$};
      \node[myrectangle] at (3.25,0) (c10) {$h^{w1}_{r1}$};
      \node[myrectangle] at (4.75,0) (c11) {$h^{w1}_{r2}$};
      \node[myrectangle] at (6.25,0) (c12) {$h^{w1}_{r3}$};
      \node[myrectangle] at (7.75,0) (c13) {$h^{w1}_{r4}$};
      \node[myrectangle] at (9.25,0) (c14) {$h^{w1}_{r5}$};
      \node[myrectangle] at (10.75,0) (c15) {$h^{w1}_{r6}$};
      \node[mytriangle] at (3.25,-1.25) (c16) {$h^{c1}_{r1}$};
      \node[mytriangle] at (4.75,-1.25) (c17) {$h^{c1}_{r2}$};
      \node[mytriangle] at (6.25,-1.25) (c18) {$h^{c1}_{r3}$};
      \node[mytriangle] at (7.75,-1.25) (c19) {$h^{c1}_{r4}$};
      \node[mytriangle] at (9.25,-1.25) (c20) {$h^{c1}_{r5}$};
      \node[mytriangle] at (10.75,-1.25) (c21) {$h^{c1}_{r6}$};

    \foreach \i/\j/\txt/\p in {
      c1/c2//below,
      c1/c3//below,
      c1/c4//below,
      c1/c5//below,
      c1/c6//below,
      c1/c7//below}
       \draw [myarrow1] (\i) -- node[sloped,font=\small,\p] {\txt} (\j);
       
    \foreach \i/\j/\txt/\p in {
      c2/c8//below,
      c3/c8//below,
      c4/c8//below,
      c5/c8//below,
      c6/c8//below,
      c7/c8//below}
       \draw [myarrow1] (\i) -- node[sloped,font=\small,\p] {\txt} (\j);

    \draw[myarrow3]
        (c2) to (c4);
    \draw[myarrow3]
        (c2.north) to[out=25,in=155] ({c5.north});
    \draw[myarrow3]
        (c2) to (c6);  
    \draw[myarrow3]
        (c2.north) to[out=25,in=155] ({c7.north}); 
    \draw[myarrow3]
        (c3) to (c4);
    \draw[myarrow3]
        (c3) to (c5);
    \draw[myarrow3]
        (c3) to (c6);
    \draw[myarrow3]
        (c3) to (c7);
    \draw[myarrow3]
        (c4) to (c6);  
    \draw[myarrow3]
        (c4.north) to[out=25,in=155] ({c7.north});
    \draw[myarrow3]
        (c5) to (c7);
    \draw[myarrow3]
        (c6) to (c7);

    \foreach \i/\j/\txt/\p in {
      c2/c9//below,
      c3/c11//below,
      c4/c12//below,
      c5/c12//below,
      c6/c13//below,
      c7/c14//below}
       \draw [myarrow2] (\i) -- node[sloped,font=\small,\p] {\txt} (\j);

    \foreach \i/\j/\txt/\p in {
      c9/c10//below,
      c10/c11//below,
      c11/c12//below,
      c12/c13//below,
      c13/c14//below,
      c14/c15//below,
      c15/c8//below,
      c16/c10//below,
      c17/c11//below,
      c18/c12//below,
      c19/c13//below,
      c20/c14//below,
      c21/c15//below,
      c16/c17//below,
      c17/c18//below,
      c18/c19//below,
      c19/c20//below,
      c20/c21//below,
      c9/c16//below,
      c10/c17//below,
      c11/c18//below,
      c12/c19//below,
      c13/c20//below,
      c14/c21//below,
      c10/c3//below,
      c11/c4//below,
      c11/c5//below,
      c12/c6//below,
      c13/c7//below}
       \draw [myarrow4] (\i) -- node[sloped,font=\small,\p] {\txt} (\j);
       
    \node[] at (1,-2.25) (n7) {};
    \node[] at (2,-2.25) (n8) {};
    \draw[myarrow4] (n7) -- (n8);
    \node[right= 0.1mm of n8] at (2,-2.25) {charging};
    
    \node[] at (4,-2.25) (n1) {};
    \node[] at (5,-2.25) (n2) {};
    \draw[myarrow2] (n1) -- (n2);
    \node[right= 0.1mm of n1] at (5,-2.25) {to charge};

    \node[] at (7,-2.25) (n9) {};
    \node[] at (8,-2.25) (n10) {};
    \draw[myarrow3] (n9) -- (n10);
    \node[right= 0.1mm of n10] at (8,-2.25) {connection};
    
    \node[] at (10,-2.25) (n5) {};
    \node[] at (11,-2.25) (n6) {};
    \draw[myarrow1] (n5) -- (n6);
    \node[right= 0.1mm of n5] at (11,-2.25) {pull-out / pull-in};

\end{tikzpicture}
\caption{Connection-based network with time-expanded charging and waiting nodes}
\label{fig:network}
\end{figure}

Our E-VSP model is defined on connection-based networks \citep{Ribeiro1994} with time-expanded charging and waiting nodes. With every depot $d \in \mathcal{D}$, we associate a network $G^d(V_d, A_d)$ with node set $V_d = \mathcal{V} \cup \{n_0^d, n_1^d\} \cup \mathcal{H}^C \cup \mathcal{H}^W$, where $n_0^d$ and $n_1^d$ represent depot $d$ at the beginning and the end of the day, respectively, and arc set $A_d$. This network contains five types of arcs $(i,j)$, namely \textit{pull-out} arcs (i.e., $(n_0^d, i)$ for $i \in \mathcal{V}$), \textit{pull-in} arcs (i.e., $(i, n_1^d)$ for $i \in \mathcal{V}$), \textit{connection} arcs (i.e.,  $(i,j)$, for $i,j \in \mathcal{V}$), \textit{charging} arcs (i.e.,  $(i,j)$, for $i \in \mathcal{H}^C \cup \mathcal{H}^W, j \in \mathcal{V} \cup \mathcal{H}^C \cup \mathcal{H}^W$), and \textit{to charge} arcs (i.e.,  $(i,j)$, for $i \in \mathcal{V}, j \in \mathcal{H}^W$). Figure \ref{fig:network} illustrates an example of such network with 6 trips $v_1$ to $v_6$, 1 charging station $h_1$, and 7 time intervals, $r_0$ to $r_6$.

Two timetabled trips $i$ and $j$ are connected by an arc $(i,j)$ in $G^d$, for all $d \in \mathcal{D}$, if and only if the departure time $d_j$ of trip $j$ is greater than or equal to the sum of the arrival time of trip $i$, the deadhead travel time $t_{ij}$ between the locations of $i$ and $j$, and a minimum layover time of $\tau$ minutes required for boarding passengers at the start of trip $j$. To minimize congestion at the terminals and reduce idle time for drivers, we implement a waiting time threshold of 45 minutes between consecutive trips $i$ and $j \in \mathcal{V}$. If the waiting time exceeds 45 minutes, the vehicle must proceed to the nearest depot after completing trip $i$. At the depot, it will remain idle until it is time to depart for the starting location of trip $j$. Since no driver attendance is required at the depot, the connection to trip $j$ remains feasible regardless of the waiting duration. We distinguish four types of charging arcs. First, for each charging station $h \in \mathcal{H}$, all consecutive nodes in $\mathcal{H}^C$ and $\mathcal{H}^W$ associated with $h$ are connected, such that two disjoint chains are formed per charging station. The nodes in $\mathcal{H}^W$ associated with the latest time interval are connected to $n_1^d$. Second, we connect each node in $\mathcal{H}^W$, except the ones associated with the last time interval, to the node in $\mathcal{H}^C$ associated with the same charging station and the next time interval. Third, an arc connects each node in $\mathcal{H}^C$ to the node in $\mathcal{H}^W$ that belongs to the same time interval and charging station. Fourth, we add an arc between a waiting node $i  \in \mathcal{H}^W$ with time interval $r + 1$ and a trip $j \in V$ if the latest departure time such that the EB arrives on time for the departure time for the trip falls in $r$ (i.e., $B_r \leq d_j - \tau - t_{ij} < E_r$). Finally, we create an arc `to charge' between a trip $i \in \mathcal{V}$ and a node $j \in \mathcal{H}^W$ if the earliest EB arrival time at the charging station falls within the corresponding time interval $r$ of the waiting node (i.e., $B_r \leq d_i + t_i +t_{ij} < E_r$).

By definition, a path in the graph $G^d$ is a sequence of pairwise compatible timetabled trips and charging activities. This path is a feasible vehicle schedule if it starts at the source node $n_0^d$, ends at the sink node $n_1^d$, and respects energy-feasibility and charging-feasibility constraints. We assume that each vehicle begins its schedule with a SoC of $\sigma^{init}$. The energy-feasibility constraints ensure that the SoC of all EBs remains above the minimum allowed SoC, denoted as $\sigma^{min}$, throughout the day. Energy consumption values, SoC values, $\sigma^{min}$, and $\sigma^{init}$ are expressed as percentages of the total battery capacity. All values are rounded to the nearest integer, resulting in admissible values within the set $\{0\%, 1\%, \dots, 100\%\}$. Moreover, charging-feasibility constraints ensure that each time an EB visits a charging station, we schedule exactly one recharging activity that can be composed of one or more consecutive nodes in $\mathcal{H}^C$. These constraints are enforced by including two charging resources in the dynamic programming algorithm, as described in Section \ref{sec:pricing}.

The cost of a vehicle schedule $s$ passing through the set of arcs $A(s)$ is given by

\begin{equation}
    c_s = \sum_{(i,j) \in A(s)} c_{ij},
\end{equation}

\noindent where $c_{ij}$ is the cost of arc $(i,j)$, consisting of a cost per vehicle used, a cost for each minute waiting outside the depot, a cost per unit of distance traveled, and a cost per charging activity. This latter cost is used, on the one hand, to discourage unnecessary visits to the charging stations when the SoC is still high and, on the other hand, to take into consideration operational costs, such as the salary of employees dedicated to each charging station.   

\subsection{Non-linear charging profile}
\label{sec:charging_profile}

The potential risks associated with excessive voltage levels, which can permanently damage the battery, can be mitigated by charging the battery using the constant current and constant voltage (CC-CV) scheme \citep{Pelletier2017}. In the CC phase, the charging current is kept constant until the voltage reaches a threshold. Then, in the CV phase, the voltage is kept constant while the current decreases. These changes in voltage and current regimes make the charging functions nonlinear. In this work, we approximate these nonlinear charging functions with piecewise linear functions, as done in \cite{Montoya2017} and \cite{Olsen2020}.

Let $\rho$ be the duration of a time interval and $x$ be the SoC at the beginning of the interval. The SoC after charging for one interval is approximated using the piecewise linear function $\lambda(x, \rho)$. If a path in $G^d$ includes $m$ consecutive charging nodes, the resulting SoC is approximated by $\lambda(x, m\rho)$---charging over a single interval of duration $m\rho$---rather than applying $\lambda(x, \rho)$ iteratively $m$ times. This approach avoids potential cumulative rounding errors. The output of $\lambda(x, \rho)$ is capped at $\sigma^{\text{max}}$, the maximum allowable SoC, or at a lower SoC bound, to prevent excessive charging at additional charging nodes.

\subsection{Chance-constrained E-VSP with battery degradation}
\label{sec:CC_model}

Consider a recommended SoC range, denoted as $[\sigma^{low}, \sigma^{up}]$, where $\sigma^{min} \leq \sigma^{low} \leq \sigma^{up} \leq \sigma^{max}$. A deterministic strategy to ensure compliance with this recommended range consists of: (1) stopping battery charging when the SoC reaches $\sigma^{up}$, and (2) scheduling EB operations such that the SoC never drops below $\sigma^{low}$, accounting for the worst-case energy consumption for all trips. While this method guarantees that the SoC remains within the prescribed range, it is conservative and may result in costly solutions. This hypothesis is examined in Section \ref{sec:results}.

Instead, we propose to limit the probability of EBs operating outside the recommended SoC range by considering the probability distribution functions (PDFs) of the energy consumption of each trip. Our goal is to equip transit agencies with a tool that ensures compliance with the recommended SoC range, typically specified by the battery supplier or leasing company, while enabling cost-efficient fleet planning. For the upper bound of the recommended SoC range, we apply a cut-off when the battery level reaches $\sigma^{up}$, similar to the deterministic approach. For the lower bound of the SoC range, we apply a chance constraint to limit the probability of falling below the recommended range. To easily refer to this approach later, we will say that a battery is \textit{overused} when its SoC is below $\sigma^{low}$.

For each trip $i \in \mathcal{V}$, we are given the probability mass function (PMF) with finite supports of the energy consumption, denoted $e_i(\mu)$. However, deterministic energy consumption is assumed for pull-out, pull-in, and deadhead trips, as these trips are typically short and do not involve passengers, thereby eliminating many sources of uncertainty in energy usage. Factors such as passenger loading, heating, and air conditioning, which significantly impact energy consumption, are not present during these trips.

Our model uses the following additional notation, which is summarized in Table \ref{tab:notation}. Let $\mathcal{S}$ be the set of all feasible vehicle schedules and $\mathcal{S}^d \subset \mathcal{S}$ be the subset of feasible vehicle schedules for depot $d \in \mathcal{D}$. Let also  $a_{is}$ be a binary parameter equal to 1 if schedule $s \in \mathcal{S}$ covers trip $i$, $y_s$ be a binary variable equal to 1 if schedule $s$ is part of the solution, and $w_s^{hr}$ be a binary parameter equal to 1 if schedule $s$ includes a charging activity at station $h$ in time interval $r$. 

\begin{table}[h!]
    \centering
    \begin{tabular}{ll}
    \toprule
       Notation  & Definition\\
    \midrule
    $\mathcal{V}$ & Set of timetabled trips\\
    $\mathcal{D}$ & Set of depots\\
    $\mathcal{H}$ & Set of charging stations\\
    $\mathcal{H}^C$ & Set of charging nodes\\
    $\mathcal{H}^W$ & Set of waiting nodes\\
    $\mathcal{R}$ & Set of time intervals\\
    $\mathcal{S}$ & Set of all feasible vehicle schedules\\
    $\mathcal{S}^d$ & Set of all feasible vehicle schedules hosted in depot $d$\\
    $y_s$ & Binary variable equal to 1 if schedule $s$ is selected\\
    $a_{is}$ & Binary parameter equal to 1 if schedule $s$ covers trip $i$\\
    $w_s^{hr}$ & Binary parameter equal to 1 if schedule $s$ includes a charging activity at station $h$ in time interval $r$\\ 
    $\tau$ & Minimum layover time \\
    $\rho$ & Duration of each time interval\\
    $\sigma^{low}$ & Lower bound on the recommended SoC range\\
    $\sigma^{up}$ & Upper bound on the recommended SoC range\\
    $\sigma^{init}$ & Initial SoC at the beginning of the planning period\\
    $\sigma^{min}$ & Minimum authorized SoC (e.g., 0\%)\\
    $\sigma^{max}$ & Maximum authorized SoC (e.g., 100\%)\\
    $b_d$ & Number of available EBs at depot $d$\\
    $g_h$ & Number of available chargers at charging station $h$\\
    $\epsilon$ & Chance constraint threshold\\
    $e_i(\mu)$ & PMF with finite supports of the energy consumption of trip $i$ (relative to the total battery capacity)\\
    \bottomrule
    \end{tabular}
    \caption{Overview of the notation of the E-VSP with battery degradation and stochastic energy consumption}
    \label{tab:notation}
\end{table}

The E-VSP with battery degradation and stochastic energy consumption can be formulated as the following chance-constrained integer program:

\begin{align}
    \label{eq:SP1}
    \text{min \qquad}&\sum_{s \in \mathcal{S}} c_s y_s \\
    \label{eq:SP2}
    \text{s.t. \qquad}&\sum_{s \in \mathcal{S}} a_{is}y_s = 1, \qquad \forall i \in \mathcal{V} \\
    \label{eq:SP3}
    &\sum_{s \in \mathcal{S}^d}  y_s \leq b_d, \qquad \forall d \in \mathcal{D}\\
    \label{eq:SP4}
    &\sum_{d \in \mathcal{D}} \sum_{s \in \mathcal{S}^d} w_{s}^{rh} y_s \leq g_h, \qquad \forall h \in \mathcal{H}, r \in \mathcal{R}\\
    \label{eq:SP5}
    &\text{Pr}\{\text{overusing the battery of one EB or more}\} \leq \epsilon\\
    \label{eq:SP6}
    &y_s \in \{0,1\}, \qquad \forall s \in \mathcal{S}.
\end{align}

The objective function (\ref{eq:SP1}) minimizes the total operational costs, while constraints (\ref{eq:SP2}) ensure that each timetabled trip is covered exactly once by a schedule, constraints (\ref{eq:SP3}) impose vehicle availability at each depot, and constraints (\ref{eq:SP4}) guarantee that for each charging station $h \in \mathcal{H}$, a maximum of $g_h$ chargers are used at the same time. Finally, the chance constraint (\ref{eq:SP5}) ensures that the probability of overusing the battery of at least one EB is less than or equal to a maximum threshold $\epsilon$.

We assume that the probabilities of battery overuse for any pair of EB schedules in $\mathcal{S}$ are independent. Let the probability $P_s = Pr\{\text{an EB remains within the safe SoC range under schedule } s\}$. Using this notation and the previous independence assumption, we can rewrite constraint (\ref{eq:SP5}) as

\begin{equation}
    1 - \prod_{s \in \mathcal{S}: y_s = 1} P_s \leq \epsilon,
\end{equation}

\noindent which, by the properties of the logs and extending the summation to all vehicle schedules, is equivalent to

\begin{equation}
    \sum_{s \in \mathcal{S}} y_s \ln(P_s) \geq \ln( 1 - \epsilon).
\end{equation}

Let $\beta_s = \ln (P_s)$ and $\beta = \ln (1 - \epsilon)$. Then, we have 

\begin{equation}
    \sum_{s \in \mathcal{S}} y_s  \beta_s \geq \beta. 
    \label{eq:SP5_modified}
\end{equation}

In what follows, we provide a method to compute $P_s$ for any given vehicle schedule $s \in \mathcal{S}$. Note that, if an EB is allowed to have a SoC below \(\sigma^{low}\), we strictly enforce for energy feasibility that the SoC remains above \(\sigma^{min}\).

\subsection{Computing the probability of overusing a battery}
\label{sec:proba_overuse}

The charging activities are interrupted when the SoC reaches $\sigma^{up}$ and EBs begin their schedules at $\sigma^{up}$. Therefore, for every EB schedule, the probability of observing a battery with a SoC greater than $\sigma^{up}$ is null. In what follows, we explain how the probability of observing a SoC smaller than $\sigma^{low}$ can be computed.

Consider a schedule $s = (0,1, \dots, n, n+1)$ where $0$ and $n+1$ are the nodes $n^d_0$ and $n^d_1$ associated with a depot $d \in \mathcal{D}$, respectively, and $1, \dots, n$ are other nodes in $V_d$. Let $X_i^s$ be the SoC at the end of node $i$ in schedule $s$ and $f_i^s(x)$ be the PMF of the event $X_i^s = x$. We are interested in $f_i^s(x)$ conditioned on the event that the EB following schedule $s$ has not yet experienced battery overuse, which we denote as $f_i^s(x|\text{not overused})$. We can compute $f_i^s(x|\text{not overused})$ recursively with $f_0^s(\sigma^{init}|\text{not overused}) = 1$ as

\begin{equation}
    \label{eq:f}
    f_i^s(x|\text{not overused}) = \begin{cases}
        \sum_{\mu = \sigma^{min}}^{\sigma^{up} - x - \iota_{i-1,i}} e_i(\mu) f_{i-1}^s(x + \mu + \iota_{i-1,i}|\text{not overused}) & \text{ if } i \in \mathcal{V} \text{ and } x \geq \sigma^{low}\\
        f_{i-1}^s(\lambda^{-1}(x,\rho)|\text{not overused}) & \text{ if } i \in \mathcal{H}^C\\
        f_{i-1}^s(x + \iota_{i-1,i}|\text{not overused}) & \text{ if } i \in \mathcal{H}^W \cup \{n_0^d\} \text{ and } x \geq \sigma^{low}\\
        0 & \text{ otherwise,}
    \end{cases}
\end{equation}

\noindent where $\lambda^{-1}(x,\rho)$ is the inverse of the charging function and outputs the initial SoC such that the final SoC after a charge of $\rho$ minutes is $x$ and $\iota_{i-1,i}$ is the deterministic energy consumption between nodes $i-1$ and $i$ due to deadheading, pulling in, or pulling out. Observe that if $\sum_{\mu = \sigma^{min}}^{\sigma^{up}} e_i(\mu) < 1$ (i.e., the energy consumption of trip $i$ exceeds the allowed battery capacity), the E-VSP with stochastic energy consumption and battery degradation becomes infeasible. Therefore, this case is excluded from consideration. The probability of staying within the recommended SoC range in schedule $s$ is thus

\begin{equation}
    P_s = \sum_{x = \sigma^{low}}^{\sigma^{up}} f_{n+1}^s(x|\text{not overused}) = F_{n+1}^s(\sigma^{up}|\text{not overused}),
    \label{eq:P_s}
\end{equation}

\noindent where $F_{i}^s(x|\text{not overused})$ is the cumulative distribution function (CDF) of $f_{i}^s(x|\text{not overused})$, computed as

\begin{equation}
\label{eq:CDF_F}
    F_i^s(x|\text{not overused}) = \sum_{x' = \sigma^{low}}^{x} f _i^s(x'|\text{not overused}), \qquad x = \sigma^{low}, \sigma^{low} + 1, \dots, \sigma^{up}.
\end{equation}

Although equation (\ref{eq:P_s}) only involves node $i = n + 1$, we define $F_i^s$ for an arbitrary node $i$ as it is used later in the labeling algorithm (see Section \ref{sec:pricing}).

\section{Heuristic branch-and-price algorithm}

\label{sec:algo}
The linear relaxation of the E-VSP formulated as a set partitioning problem provides a tight lower bound but involves a large number of decision variables, one for each EB schedule. Therefore, enumerating all possible schedules is neither practical nor efficient. To address this, a column generation algorithm \citep[see, e.g.,][]{CC2005, Lubbecke2005} is integrated into the branch-and-bound algorithm to dynamically generate variables as needed. This combined approach is known as branch-and-price \citep{Barnhart1998, Costa2019}.

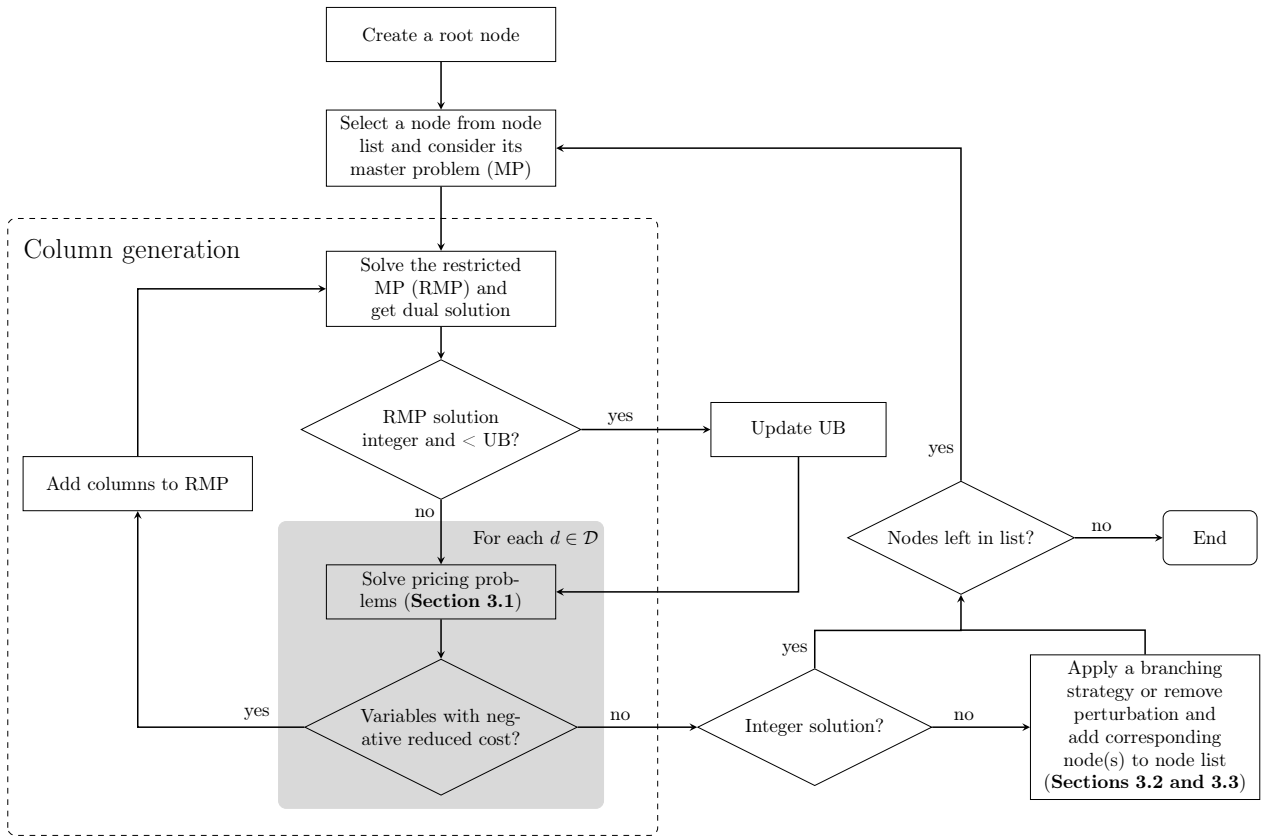
\begin{figure}[!ht]
\centering
\scalebox{0.72}{
    \begin{tikzpicture}[node distance=2.1cm]

    \draw[rounded corners, fill=black!15, black!15] (-3,-14.3) rectangle ++(6,5.3);
    \node[draw=none] at (1.75,-9.3) {For each $d \in \mathcal{D}$};
    \draw[dashed, rounded corners, draw=black] (-8,-14.8) rectangle ++(12,11.4);
    \node[draw=none] at (-5.7,-4.0) {\Large Column generation};

    \node (pro2) [process] {Create a root node};
    \node (pro3) [process, below of=pro2] {Select a node from node list and consider its master problem (MP)};
    \node (pro4) [process, below of=pro3, yshift = -0.5cm] {Solve the restricted MP (RMP) and get dual solution};
    \node (pro5) [decision, below of=pro4, yshift=-0.5cm] {RMP solution integer and < UB?};
    \node (pro6) [process, below of=pro5, yshift=-0.9cm] {Solve pricing problems (\textbf{Section \ref{sec:pricing}})};
    \node (dec7) [decision, below of=pro6, yshift=-0.4cm] {Variables with negative reduced cost?};
    \node (pro8a) [process, left of=pro5, xshift=-3.5cm, yshift=-1cm] {Add columns to RMP};
    \node (dec8b) [decision, right of=dec7, xshift=4.8cm] {Integer solution?};
    \node (pro9b) [smallprocess, right of=pro5, xshift=4.5cm] {Update UB};
    \node (pro10b) [process, right of=dec8b, xshift=4cm] {Apply a branching strategy or remove perturbation and add corresponding node(s) to node list (\textbf{Sections \ref{sec:branching} and \ref{sec:perturbation}})};
    \node (dec11) [decision, right of=pro6, xshift=7.5cm, yshift=1cm] {Nodes left in list?};
    \node (end) [startstop, right of=dec11, xshift=2.5cm] {End};

    \draw [arrow] (pro2) -- (pro3);
    \draw [arrow] (pro3) -- (pro4);
    \draw [arrow] (pro4) -- (pro5);
    \draw [arrow] (pro5) -- node[anchor=south, xshift=-0.3cm, yshift=+0.15cm] {no} (pro6);
    \draw [arrow] (pro6) -- (dec7);
    \draw [arrow] (dec7) -- node[anchor=east, yshift=+0.2cm, xshift=+0.0cm] {no} (dec8b);
    \draw [arrow] (dec7) -| node[anchor=south, xshift=2.2cm] {yes} (pro8a);
    \draw [arrow] (pro8a) |- (pro4);
    \draw [arrow] (pro5) -- node[anchor=east, yshift=+0.2cm, xshift=-0.1cm] {yes} (pro9b);
    \draw [arrow] (dec8b)-- node[anchor=east, yshift=+0.2cm, xshift=0.0cm] {no} (pro10b);
    \draw [arrow] (dec11) -- node[anchor=east, yshift=+0.2cm, xshift=+0cm] {no} (end);
    \draw[arrow] (dec8b.north) -- node[anchor=east, yshift=+0cm, xshift=0cm] {yes}  ++(0,0.7)  -| ++(2.7,0) -- (dec11.south);
    \draw[thick] (pro10b.north) -- ++(0,0.45) -| ++(-3.4,0);
    \draw [arrow] (dec11) |- node[anchor=east, xshift=0cm, yshift=-5.6cm] {yes} (pro3);
    \draw [arrow] (pro9b) |- (pro6);

    \end{tikzpicture}
}
\caption{Flowchart of the branch-and-price algorithm \citep[from][]{Ricard2024}}
\label{fig:branch_and_price}
\end{figure}

The framework of our stochastic branch-and-price algorithm is illustrated in Figure \ref{fig:branch_and_price}. The branch-and-bound framework handles the initialization of the root node (i.e., the linear relaxation of (\ref{eq:SP1})-(\ref{eq:SP6}), hereafter referred to as the \textit{master problem} (MP)), the management of the list of active nodes, the update of the upper bound (UB), and the handling of the branching strategy and perturbation method. The UB is updated when an integer solution to the restrcited MP (RMP)---the MP limited to a subset \( \mathcal{S}' \subseteq \mathcal{S} \) of the schedule variables \( y_s \)---with a cost less than the current UB is found. Each node in the branch-and-bound tree is solved using a tailored column generation algorithm.

In each iteration of the column generation algorithm, a RMP is solved. Given the dual solution to the RMP, a pricing problem is then solved to generate new columns. For the E-VSP with battery degradation and stochastic energy consumption, the pricing problem is separable by depot and consists in solving a shortest path problem with stochasticity \citep{Boland2015, Wellman2013} in \( G^d \) for each \( d \in \mathcal{D} \), where arc costs are adjusted based on the most recent dual solution of the RMP. If variables with negative reduced costs are identified, they are added to \( \mathcal{S}' \), triggering a new iteration of the column generation algorithm. If no such variables are found, the column generation algorithm terminates, and the current RMP solution is guaranteed to be optimal for the MP of the corresponding branch-and-bound node.

Details of the stochastic dynamic programming algorithm used to solve the pricing problems, including label components, extension functions, and dominance rules, are provided in Section \ref{sec:pricing}. The branching strategy and perturbation method employed to enhance computational efficiency are discussed in Sections \ref{sec:branching} and \ref{sec:perturbation}, respectively.

\subsection{Stochastic pricing problems}
\label{sec:pricing}

Let $(u_i)_{i \in \mathcal{V}}$, $(\pi_d)_{d \in \mathcal{D}}$, $(\alpha^{hr})_{h \in \mathcal{H}, r \in \mathcal{R}}$, and $\theta$ be the dual variables associated with constraints (\ref{eq:SP2}), (\ref{eq:SP3}), (\ref{eq:SP4}), and (\ref{eq:SP5}), respectively. Each shortest path problem with stochasticity is defined on a network $G^d$ with the cost of each arc modified as presented in Table \ref{table:cost_of_arcs}.

\begin{table}[!t]
    \centering
    \begin{tabular}{llll}
    \toprule
        Name & Start node $i$ & End node $j$ & $\Tilde{c}^s_{ij}$\\
    \midrule
        Pull-out & $i = n_o^d$ & $ j \in \mathcal{V}$ & $c^s_{ij} - \pi_d$\\
        Connection & $ i \in \mathcal{V}$ & $j \in \mathcal{V} \cup \mathcal{H}^W$  & $c^s_{ij} - u_i$\\
        Charging (from charg.)& $i = h_r^c \in \mathcal{H}^C$ & $j \in \mathcal{H}^C \cup \mathcal{H}^W$ & $c^s_{ij} - \alpha^{hr}$\\
        Charging (from wait.) & $i \in \mathcal{H}^W$ &  $j \in V \cup \mathcal{H}^C \cup \mathcal{H}^W $& $c^s_{ij}$ \\
        To charge & $ i \in \mathcal{V}$ & $ j \in \mathcal{H}^W$ & $c^s_{ij}$ \\
        Pull-in (from terminal) & $i \in \mathcal{V}$ & $j = n_1^d$ & $c^s_{ij} - u_i - \theta \> \ln P_s$\\
        Pull-in (from wait.)  & $i \in \mathcal{H}^W$ & $j = n_1^d$ & $c^s_{ij} - \theta \> \ln P_s$\\
    \bottomrule
    \end{tabular}
    \caption{Modified arc costs}
    \label{table:cost_of_arcs}
\end{table}

A feasible schedule is a path in $G^d$ starting at node $n_0^d$ and ending at node $n_1^d$. The reduced cost $\Tilde{c}_s$ of schedule $s$ is therefore given by

\begin{align}
    \Tilde{c}_s & = c_s - \sum_{i \in \mathcal{V}} a_{is} u_i - \pi_d - \sum_{h \in \mathcal{H}} \sum_{r \in \mathcal{R}} w_s^{hr} \alpha^{hr} - \theta \> \ln P_s.
\end{align}

The component of the reduced cost associated with the chance constraint (\ref{eq:SP5}), i.e., $\theta \> \ln P_s$, is not decomposable by arc. This is because the probability of overusing the battery in schedule $s$ is only revealed once the entire schedule is constructed. Hence, we add this cost to the last arc of a schedule, i.e., the pull-in arc.

Each path is constructed using a dynamic programming approach, specifically a stochastic labeling algorithm \citep{Wellman2013, Boland2015, Ricard2024}. The algorithm operates as follows: starting from an initial label at the source node \(n^d_0\), labels are propagated through the network using resource extension functions to construct partial paths. This process continues until the sink node is reached, where the paths represent complete schedules. Each label encodes the necessary information to evaluate the feasibility of any of its extension and to compute the contribution of the chance constraint to the reduced cost on the pull-in arc at the end of the schedule. Infeasible paths are discarded as soon as they are identified as such. To improve efficiency, a dominance procedure is applied to compare labels, eliminating non-Pareto-optimal labels. The main components of the stochastic labeling algorithm used in our tailored branch-and-price algorithm, namely the labeling procedure and the dominance rule, are detailed next.

\subsubsection{Labeling procedure}

Each label encodes a representation of the probability of overusing the battery of the EB, the worst-case SoC, the accumulated reduced cost, and two additional resources. To compute the probability of overusing the battery of an EB along a path $p$ ending at node $i$ and its extensions, it is necessary to store $F_i^p(x|\text{not overused})$ for integer values of $x$ between $\sigma^{low}$ and $\sigma^{up}$. Simply relying on the CDF of the SoC is insufficient because it lacks memory. By definition, $F_i^p(x|\text{not overused})$ incorporates this memory by conditioning the CDF of the SoC by the event that the EB following path $p$ has not yet experienced battery overuse. The two additional resources are for planning feasible charging activities. The first, $R_i^p$, represents the number of completed charging activities since the last timetabled trip in $p$. It ensures that at most one recharge occurs for each visit at a charging station. The second, $E_i^p$, tracks the number of times an EB has waited at a charging station without being recharged since the last timetabled trip in $p$. This resource enforces that a recharge occurs during such a visit. Both resources are initialized to 0.

A label is represented as \( L_i^{p'} = (F_i^{p'}(\sigma^{low}|\text{not overused}), \dots, F_i^{p'}(\sigma^{up}|\text{not overused}), \omega_i^p, C_i^{p'}, R_i^{p'}, E_i^{p'}) \), where $\omega_i^{p'}$ is the worst-case SoC for path $p'$, which starts at the origin node $n_o^d$ and ends at node $i$, and $C_i^{p'}$ is the accumulated reduced cost for path $p'$. The extension of $L_i^{p'}$ along the arc $(i,j) \in A_d$ is done by applying the following extension functions to create a new label $L_j^p$ at node $j$.

Let $\sigma_i^p$ be the SoC at node $i$ of path $p$. Assuming each EB starts the day fully recharged (i.e., $\sigma_i^p = \sigma^{up}$ for $i = n_o^d$), the extension function to compute $F_j^p(x|\text{not overused})$ based on $F_i^p(x|\text{not overused})$, for $x = \sigma^{low}, \sigma^{low} + 1, \dots, \sigma^{up}$, is derived from equations (\ref{eq:f}) and (\ref{eq:CDF_F}), and given by

\begin{equation}
\label{eq:ext_F}
    F_j^p(x|\text{not overused}) = \begin{cases}
        \sum_{\mu =\sigma^{min}}^{\sigma^{up}- x - \iota_{i-1,i}} e_j(\mu) F_i^{p'}(x + \mu + \iota_{i,j}|\text{not overused}) & \text{ if } j \in \mathcal{V}\\
        F_i^{p'}(\lambda^{-1}(x,\rho)) & \text{ if } j \in \mathcal{H}^C\\
        F_i^{p'}(x + \iota_{i,j}) & \text{ if } i \in \mathcal{H}^W \cup \{n_0^d\}\\
        0 & \text{ otherwise.} 
    \end{cases}
\end{equation}

Furthermore, $\omega_j^{p}$ can be computed using

\begin{equation}
\label{eq:ext_xi}
\omega_j^{p} = \begin{cases}
        \omega_i^{p'} - \max \{\mu | e_j(\mu)>0\} - \iota_{i,j} & \text{ if } j \in \mathcal{V}\\
        \lambda(\omega_i^{p'}, \rho)& \text{ if } j \in \mathcal{H}^C\\
        \omega_i^{p'} - \iota_{i,j}  & \text{ otherwise,}
    \end{cases}
\end{equation}

\noindent where $\max \{\mu | e_j(\mu)>0\}$ outputs the largest energy consumption value with positive probability. The accumulated reduced cost $C_j^p$ is updated as

\begin{equation}
\label{eq:ext_cost}
    C_j^p = C_i^{p'} + \Tilde{c}_{ij}^p,
\end{equation}

\noindent and, lastly, the resources $R_{i}^{p'}$ and $E_{i}^{p'}$ are extended as

\begin{equation}
\label{eq:ext_R}
    R_{j}^{p} = \begin{cases}
        R_{i}^{p'} + 1 & \text{ if } i \in \mathcal{H}^W \text{ and } j  \in \mathcal{H}^C\\
        R_{i}^{p'} - 1  & \text{ if } i  \in \mathcal{H}^W \text{ and } j \in \mathcal{V}\\
        R_{i}^{p'} & \text{ otherwise,}
    \end{cases}
\end{equation}

\noindent and

\begin{equation}
\label{eq:ext_T}
    E_{j}^{p} = \begin{cases}
        E_{i}^{p'} + 1 & \text{ if } i \in \mathcal{V} \text{ and } j \in \mathcal{H}^W\\
        E_{i}^{p'} - 1 & \text{ if } i \in \mathcal{H}^W \text{ and } j \in \mathcal{H}^C\\
        E_{i}^{p'} & \text{ otherwise.}
    \end{cases}
\end{equation}

To guarantee energy feasibility, a path \( p \) ending at node \( i \) is discarded if there is a positive probability of running out of energy, specifically if \( \omega_i^p < \sigma^{min} \). Additionally, any path for which the chance constraint (\ref{eq:SP5}) would be automatically violated—--whether for the path itself or any of its potential extensions---is excluded. This occurs when \( \ln(F_i^p(\sigma^{up}|\text{not overused})) < \ln(1 - \epsilon) \). Furthermore, a path \( p \) ending at node \( i \) is discarded if the resource values \( E_i^p \) or \( R_i^p \) exceed their respective resource windows. These windows require $E_i^p$ and $R_i^p$ to be less than or equal to $0$ if $i \in \mathcal{V}$ and less than or equal to $1$ otherwise.

\subsubsection{Label dominance}

Consider two paths $p_1$ and $p_2$, both ending at node $i$. Path $p_1$ dominates path $p_2$, and thus $p_2$ can be discarded, when the following conditions hold:

\begin{enumerate}[label=(\roman*)]
    \item $C_i^{p_1} \leq C_i^{p_2}$
    \item $R_i^{p_1} \leq R_i^{p_2}$
    \item $E_i^{p_1} \leq E_i^{p_2}$
    \item $\omega_i^{p_1} \geq \omega_i^{p_2}$
    \item $ F_{i}^{p_1}(\sigma^{up} |\text{not overused}) \geq F_{i}^{p_2}(\sigma^{up}|\text{not overused})$
    \item $ F_{i}^{p_1}(\sigma^{up}|\text{not overused}) - F_{i}^{p_1}(x - 1|\text{not overused}) \geq F_{i}^{p_2}(\sigma^{up}|\text{not overused}) - F_{i}^{p_2}(x -1|\text{not overused})$, for all $x \in \{\sigma^{low}, \sigma^{low} + 1, \dots, \sigma^{up} - 1\}$
\end{enumerate}

The first four conditions are standard in shortest path problems with resource constraints. Conditions (v) and (vi), however, are dominance rules tailored to our problem. On the one hand, condition (v) compares the probabilities that the SoC stays within the recommended range in $p_1$ and $p_2$. It ensures that paths with a higher probability of maintaining the SoC within this range are preferred. Condition (vi), on the other hand, evaluates the complementary CDFs of the SoC conditioned on the event that the EBs following paths $p_1$ and $p_2$ remain within the recommended SoC range to prioritize paths with higher SoC values. Higher SoC values are desirable as they reduce the risk of falling below \(\sigma^{low}\) as paths are extended. Note that although conditions (v) and (vi) are based on probability distributions, they remain exact because we propagate the full distributions and apply the dominance rule considering all possible realizations.

\subsection{Branching strategy}
\label{sec:branching}

To limit the time spent exploring the branch-and-bound nodes and to identify integer solutions within a reasonable timeframe, we implement a heuristic branching strategy, namely a diving heuristic. First, we branch on the total number of vehicles used by creating one or two child nodes with upper and lower bounds on the depot capacity. Second, we apply three rounding strategies: (i) rounding schedule variables, (ii) rounding \textit{connection} arcs, and (iii) rounding any arc (including \textit{to charge} and \textit{charging} arcs).

When applying strategy (i), a single branch-and-bound node is created by fixing one or more schedule variables to 1. Variables with the largest fractional values are prioritized for rounding. At each node, a maximum of three variables with fractional parts greater than or equal to 0.99 are selected. If fewer than three variables meet this threshold,
all eligible variables are selected. If no variables exceed this threshold, the variable with the largest fractional part (below 0.99) is chosen. The threshold and the number of variables selected were determined empirically during a preliminary experimental campaign.

The same procedure is used for strategies (ii) and (iii), except that arcs are fixed instead of variables. Fixing an arc $(i,j)$ to 1 is achieved by removing all arcs \((i,k)\) and \((k,j) \in A_d\), \(d \in \mathcal{D}\), such that \(k \in V_d \setminus \{i,j\}\).

The rounding strategies (i)--(iii) are alternated from one node to the next after the initial branching on the total number of vehicles, with priority given to strategies (i) and (ii). At each node, the strategy with the highest score is selected. The score is defined as the largest fractional value among the selected variables or arcs in strategies (i) and (ii). For strategy (iii), the score is defined as the largest fractional value among the selected arcs, multiplied by 0.7.

\subsection{Constraint perturbation}
\label{sec:perturbation} 

E-VSPs formulated as a set partitioning problem are highly degenerate, and this degeneracy generally increases with the time horizon considered or, equivalently, the average number of trips per schedule \citep[see ][]{Oukil2007,Benchimol2012}. We use the simple constraint perturbation strategy of \cite{Charnes1952} to reduce the degeneracy in our experiments. Let $\eta_i^+$ and $\eta_i^-$ be perturbation variables that allow an under- and over-covering of trip $i \in \mathcal{V}$ of up to $\xi_i^+$ and $\xi_i^-$, respectively. The perturbed MP is given by

\begin{align}
    \label{eq:SP1_pert}
    \text{min \qquad}&\sum_{s \in \mathcal{S}} c_s y_s + \sum_{i \in \mathcal{V}} (\delta^{+}_i \eta^+_i +\delta^{-}_i \eta^-_i) \\
    \label{eq:SP2_pert}
    \text{s.t. \qquad}&\sum_{s \in \mathcal{S}} a_{is}y_s + \eta_i^+ - \eta_i^- = 1, \qquad \forall i \in \mathcal{V} \\
    \label{eq:SP3_pert}
    &(\ref{eq:SP3}), (\ref{eq:SP4}, (\ref{eq:SP5_modified}),\nonumber \\
    &0 \leq \eta_i^+ \leq \xi_i^+, \qquad \forall i \in \mathcal{V}\\
    \label{eq:SP7_pert}
    &0 \leq \eta_i^- \leq \xi_i^-, \qquad \forall i \in \mathcal{V}\\
    \label{eq:SP8_pert}
    &0 \leq y_s \leq 1, \qquad \forall s \in \mathcal{S},
\end{align}

\noindent where $\delta_i^+$ and $\delta_i^-$ are the penalties for under- and over-covering trip $i \in \mathcal{V}$, respectively.

\section{Experimental results}
\label{sec:results}

We evaluate our approach using simulated instances of one EB line in the city of Montr\'eal, which consists of 42 stops covering a distance of 8.5 kilometers. The random instances have approximately 60, 150, and 250 trips, each operating in two possible directions and evenly distributed over a 19-hour period from 5:00 AM to midnight. 

We conduct five tests for each instance size, denoted as I1, I2, and I3, where the start times of the trips are slightly shifted in each of the five tests. Table \ref{table:instances} presents the average characteristics of the families of instances I1, I2, and I3, namely the instance family name (Instance), the average number of timetabled trips ($|\mathcal{V}|$), the average number of arcs ($|A|$), and the number of chargers per charging station (Char. station capacity). In our tests, we address the single-depot and single-charging station case (i.e., $|\mathcal{D}| = 1$ and $|\mathcal{H}| = 1$), although the model and algorithm presented can easily be applied to cases where the number of depots or charging stations is greater than or equal to $2$.

\begin{table}[ht!]
\caption {Average properties of the families of instances I1 - I3}
\label{table:instances}
\begin{center}
\begin{tabular}[t]{l rrrrr}
    \toprule
    Instance&$|\mathcal{V}|$&$|A|$&Char. station\\
    &&&capacity\\
    \midrule
    I1&60&2,190&1\\ 
    I2&155&11,900&2\\ 
    I3&248&29,499&3\\ 
    \bottomrule
\end{tabular}
\end{center}
\end{table}

We assume that the energy consumption rate (kWh/km) follows a normal distribution, consistent with findings from previous studies \citep[see][]{Bie2021, Abdelaty2021}. For each trip in $\mathcal{V}$, the mean of the energy consumption rate distribution is sampled from an exponential distribution $\text{Exp}(1.57, 0.26)$, where 1.57 and 0.26 represent the location and scale parameters, respectively. This distribution has a mean value of 1.83 kWh/km. The variance of the normal distribution is determined by sampling a value from a uniform distribution $\text{U}(0.35, 0.5)$. The parameters of the exponential and uniform distributions were calibrated using data from \cite{Basma2020}.

The following model parameters are used in our tests. The costs per vehicle used, per minute of travel (excluding the travel time of timetabled trips), per minute of waiting outside the depot, and per charging activity are set to 1,000, 0.4, 0.2, and 10, respectively. The penalties for under- and over-covering trip $i \in \mathcal{V}$, $\delta_i^+$ and $\delta_i^-$, are set to $1$ for every trip $i \in \mathcal{V}$ and the upper bounds $\xi_i^+$ and $\xi_i^-$ are randomly sampled from the uniform distribution U(0,0.1) for every trip. We consider 12-meter single-deck EBs with 300 kWh batteries, the same type of EBs as those studied in \cite{Basma2020} to ensure consistency, and time intervals of $\rho = 15$ minutes. Fast chargers with an approximate power output of 450 kW are used. These chargers deliver 7.5 kWh per minute for SoC levels between 0\% and 80\%, 6 kWh per minute for SoC levels between 80\% and 90\%, and 3.75 kWh per minute for SoC levels between 90\% and 100\%. Consequently, the chargers can fully charge a 300 kWh battery in 45 minutes.

We conduct our experiments on a Linux machine equipped with 16 Intel Xeon ES-2637 v4 processors running at 3.50 GHz and a RAM of 125 GB. The branch-and-price algorithm is implemented using the GENCOL library, version 4.5, and all linear programs are solved by the commercial solver CPLEX 22.1.

The remainder of this section is organized as follows. In Section \ref{sec:performance_det}, we analyze the algorithm behavior on I1, I2, and I3 for the deterministic and stochastic cases. Section \ref{sec:comparison_det_and_sto} compares our chance-constrained model for the E-VSP with battery degradation and stochastic energy consumption to the deterministic approach in terms of operational costs and number of EBs needed. Section \ref{sec:selection_epsilon} presents a potential method for selecting the probability threshold ($\epsilon$) for the chance constraint. It is worth noting that, since our approach does not rely on a specific capacity fading model, alternative methods can be employed to choose an $\epsilon$ value that aligns with the unique constraints and objectives of each public transportation agency, as well as to comply with the guidelines provided by the battery leasing company or the manufacturer.

\subsection{Heuristic performance}
\label{sec:performance_det}

For a recommended SoC range of $[\sigma^{low}, \sigma^{up}]$, we consider a deterministic approach in which SoC values below $\sigma^{\text{low}}$ are strictly prohibited. In our framework, this corresponds to setting $\sigma^{\text{low}} = \sigma^{\text{min}}$, or equivalently, $\epsilon = 0\%$. In this case, the chance constraint (\ref{eq:SP5_modified}) becomes redundant. Energy consumption PDFs, which are used to compute the probability of battery overuse, can be replaced by a single value representing the worst-case energy consumption for each trip. This deterministic approach guarantees energy feasibility under the most conservative assumptions and serves as a baseline for comparison in the subsequent analysis.

Table \ref{table:heur_perform} summarizes the heuristic performance of our algorithm. The left-hand side (`E-VSP') presents results for the deterministic case, while the right-hand side (`Stochastic E-VSP') displays results for the E-VSP with battery degradation and stochastic energy consumption. The table reports the average relative difference in percentage between the UB and the lower bound (Gap), the number of branching nodes explored (BBn), and the computing times (CPU time) in seconds, including the total CPU time (Total), the time to solve the root node (Root), and the time to solve the pricing problems (Pricing). For the deterministic case, these metrics are averaged over five tests per instance family and two recommended SoC ranges, namely 20–80\% and 30–80\%. For the stochastic case, the metrics are averaged over five tests per instance family, 11 $\epsilon$ values, where \(\epsilon \in \{ 0.001, 0.005, 0.01, 0.05, 0.10, 0.15, 0.20, 0.25, 0.30, 0.40, 0.50 \}\), and the two same recommended SoC ranges. This amounts to a total of 110 tests per instance family for the stochastic case. 

Our branch-and-price heuristic achieves near-optimal solutions for both the deterministic and stochastic cases of the E-VSP, with average optimality gaps of 0.03\% and 0.04\%, respectively. Specifically, the gaps remain small, with the highest observed values being 0.04\% for the deterministic case and 0.06\% for the stochastic case. The majority of the CPU time is spent solving the pricing problems. In the deterministic E-VSP, between 74.5\% and 79.9\% of the total CPU time is spent on pricing on average, while in the stochastic E-VSP, this proportion increases to between 90.9\% and 96.6\%, depending on the instance family. The time spent solving the root node varies significantly across instance families, ranging from 5.3\% to 54.0\% of the total CPU time on average. In the deterministic case, instance I3 exhibits the lowest proportion of time allocated to solving the root node on average, likely because it also has the highest average number of nodes explored, at 2,633. Additionally, we observe that the total CPU time increases significantly with the instance size. For instance, the time grows from 3.2 seconds (I1) to 4,424.7 seconds (I3) in the deterministic case and from 11.5 seconds (I1) to 23,911.5 seconds (I3) in the stochastic case. The growth rate is more pronounced for the E-VSP with battery degradation and stochastic energy consumption compared to the deterministic baseline, highlighting the added complexity introduced by uncertainty in the stochastic model.

\begin{table}[!ht]
\caption {Average heuristic performance of the deterministic E-VSP and the E-VSP with battery degradation and stochastic energy consumption}
\label{table:heur_perform}
\begin{center}
\begin{tabular}[t]{l llllll llllll}
    \toprule
    &\multicolumn{5}{c}{E-VSP}&\multicolumn{5}{c}{Stochastic E-VSP}\\
    \cmidrule(lr){2-6}\cmidrule(lr){7-11}
    &&&\multicolumn{3}{c}{CPU times (s)}&&&\multicolumn{3}{c}{CPU times (s)}\\
    \cmidrule(lr){4-6}\cmidrule(lr){9-11}
    Inst.&Gap&BBn&Total&Root&Pric.&Gap&BBn&Total&Root&Pric.\\
    &(\%)&&&(\%)&(\%)&(\%)&&&(\%)&(\%)\\
    \midrule
    I1&0.02&52&3.2&54.0&79.9&0.02&73&11.5&59.3&90.9\\
    I2&0.03&148&92.6&37.3&75.9&0.05&443&1,703.8&50.0&96.4\\
    I3&0.04&2,633&4,424.7&5.3&74.5&0.06&1,314&23,911.5&26.4&96.6\\
    \bottomrule
\end{tabular}
\end{center}
\end{table}

\subsection{Comparison of the deterministic and stochastic E-VSP with battery degradation}
\label{sec:comparison_det_and_sto}

This section compares a deterministic approach to the E-VSP to our stochastic approach considering battery degradation. Table \ref{table:results} illustrates the trade-off between the operational costs and the probability of not complying with the recommended SoC range in one vehicle schedule or more. Each value in the table represents the average over five tests conducted for each instance family. The columns include the instance family name (Inst.), the chance constraint threshold ($\epsilon$), and, for the recommended SoC ranges of 20-80\% and 30-80\%, the corresponding operational costs (Op. costs), the improvement in operational costs relative to the deterministic baseline (Op. costs impro.), and the number of EBs used (\# EBs). For each instance family, the first row (with $\epsilon = 0$) corresponds to the deterministic benchmark.

\begin{table}[!ht]
\caption {Average operational costs vs. probability of not complying with the recommended SoC range in one vehicle schedule or more}
\label{table:results}
\begin{center}
\begin{tabular}[t]{ll lll lll}
    \toprule
    &&\multicolumn{3}{c}{Range 20\% - 80\%}&\multicolumn{3}{c}{Range 30\% - 80\%}\\
    \cmidrule(lr){3-5}\cmidrule(lr){6-8}
    Inst.&$\epsilon$&Op. costs&Op. costs&\# EBs& Op. costs&Op. costs&\# EBs\\
    &&&impro. (\%)&&&impro. (\%)&\\
    \midrule
    I1&0&4,102.2&0.00&3.8&4,312.7&0.00&4.0\\
    &0.001&3,896.8&5.01&3.6&4,103.1&4.86&3.8\\
    &0.005&3,896.5&5.01&3.6&4,100.8&4.91&3.8\\
    &0.01&3,896.2&5.02&3.6&4,098.3&4.97&3.8\\
    &0.05&3,895.2&5.05&3.6&3,898.0&9.62&3.6\\
    &0.10&3,894.3&5.07&3.6&3,897.5&9.63&3.6\\
    &0.15&3,894.3&5.07&3.6&3,897.2&9.63&3.6\\
    &0.20&3,894.3&5.07&3.6&3,897.3&9.63&3.6\\
    &0.25&3,894.3&5.07&3.6&3,897.1&9.64&3.6\\
    &0.30&3,894.3&5.07&3.6&3,897.0&9.64&3.6\\
    &0.40&3,894.3&5.07&3.6&3,896.6&9.65&3.6\\
    &0.50&3,894.3&5.07&3.6&3,896.4&9.65&3.6\\[0.2cm]

    I2&0&8,313.6&0.00&7.8&8,540.2&0.00&8.0\\
    &0.001&8,275.6&0.46&7.8&8,508.0&0.38&8.0\\
    &0.005&8,271.9&0.50&7.8&8,310.2&2.69&7.8\\
    &0.01&8,271.4&0.51&7.8&8,305.3&2.75&7.8\\
    &0.05&8,066.0&2.98&7.6&8,293.6&2.89&7.8\\
    &0.10&8,065.6&2.98&7.6&8,288.5&2.95&7.8\\
    &0.15&8,063.5&3.01&7.6&8,287.0&2.97&7.8\\
    &0.20&8,064.8&2.99&7.6&8,282.9&3.01&7.8\\
    &0.25&8,063.5&3.01&7.6&8,283.9&3.00&7.8\\
    &0.30&8,063.5&3.01&7.6&8,279.2&3.06&7.8\\
    &0.40&8,064.4&3.00&7.6&8,276.4&3.09&7.8\\
    &0.50&8,064.4&3.00&7.6&8,274.4&3.11&7.8\\[0.2cm]
    
    I3&0&12,157.7&0.00&11.4&12,782.3&0.00&12.0\\
    &0.001&11,889.6&2.21&11.2&12,352.8&3.36&11.6\\
    &0.005&11,691.6&3.83&11.0&11,971.1&6.35&11.2\\
    &0.01&11,687.7&3.87&11.0&11,958.4&6.45&11.2\\
    &0.05&11,676.3&3.96&11.0&11,932.1&6.65&11.2\\
    &0.10&11,673.4&3.98&11.0&11,923.4&6.72&11.2\\
    &0.15&11,673.9&3.98&11.0&11,915.0&6.79&11.2\\
    &0.20&11,673.7&3.98&11.0&11,907.8&6.84&11.2\\
    &0.25&11,672.2&3.99&11.0&11,903.1&6.88&11.2\\
    &0.30&11,673.8&3.98&11.0&11,901.1&6.89&11.2\\
    &0.40&11,674.0&3.98&11.0&11,898.6&6.91&11.2\\
    &0.50&11,674.0&3.98&11.0&11,710.0&8.39&11.0\\
    \bottomrule
\end{tabular}
\end{center}
\end{table}

For all instance families (I1, I2, I3), even a small increase in \(\epsilon\) leads to important reductions in operational costs. This reduction is primarily attributed to the decrease in the number of EBs required, which is a desired outcome. For example, in I3, increasing $\epsilon$ from \(\epsilon = 0\) to \(\epsilon = 0.005\) reduces the number of EBs required from 12.0 to 11.2, resulting in a 6.35\% improvement in operational costs for the 30–80\% range. Similar trends are observed in I1 and I2, where operational costs almost always decrease monotonically as the chance constraint is slightly relaxed.

\begin{figure}[ht!]
    \centering
    \subfloat[Range 20-80\%]{\label{fig:results_I1_20-80}\includegraphics[width=0.48\textwidth, height = 6.5 cm, keepaspectratio]{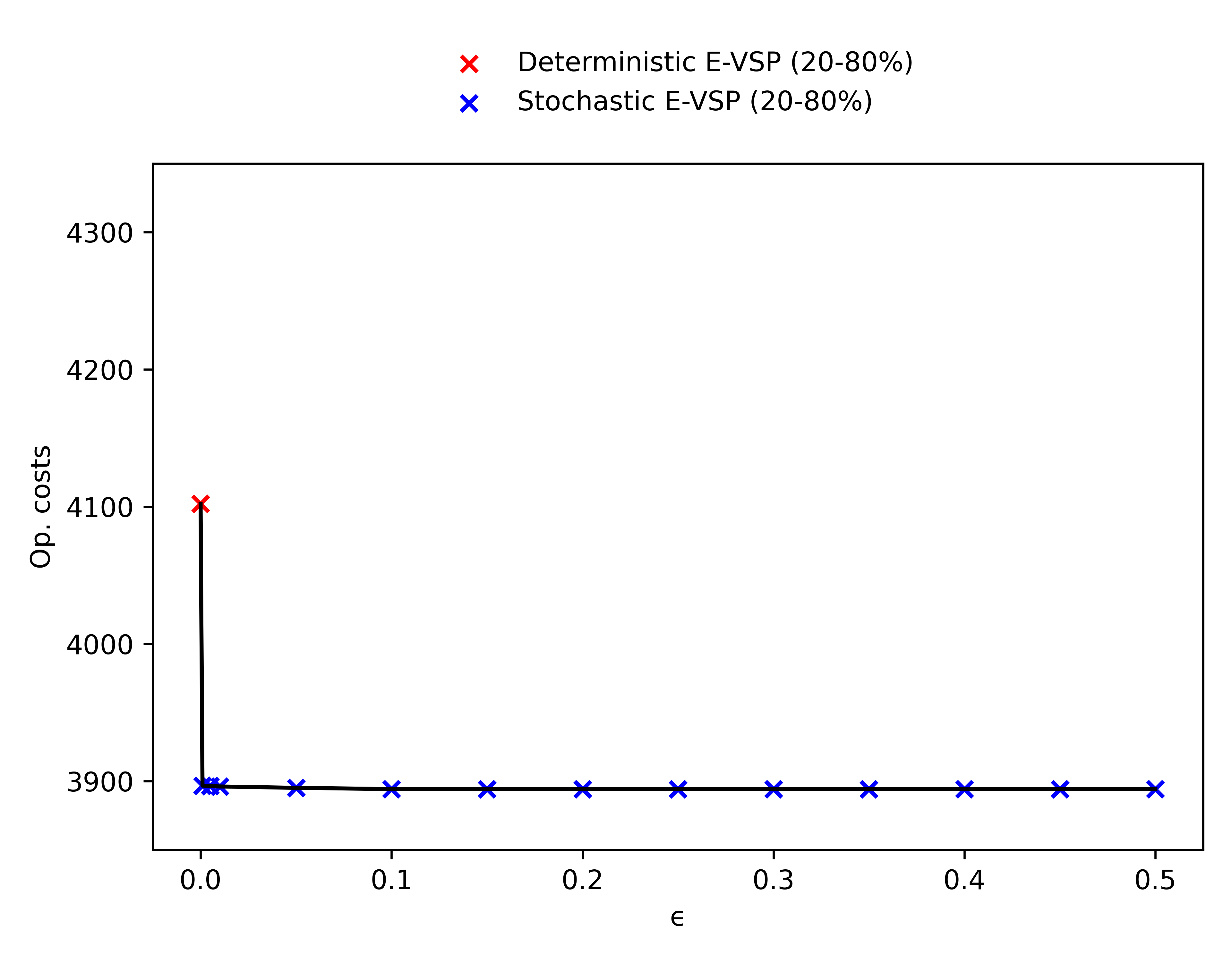}}\hspace{0cm}
    \subfloat[Range 30-80\%]{\label{fig:results_I1_30-30}\includegraphics[width=0.48\textwidth, height = 6.5 cm, keepaspectratio]{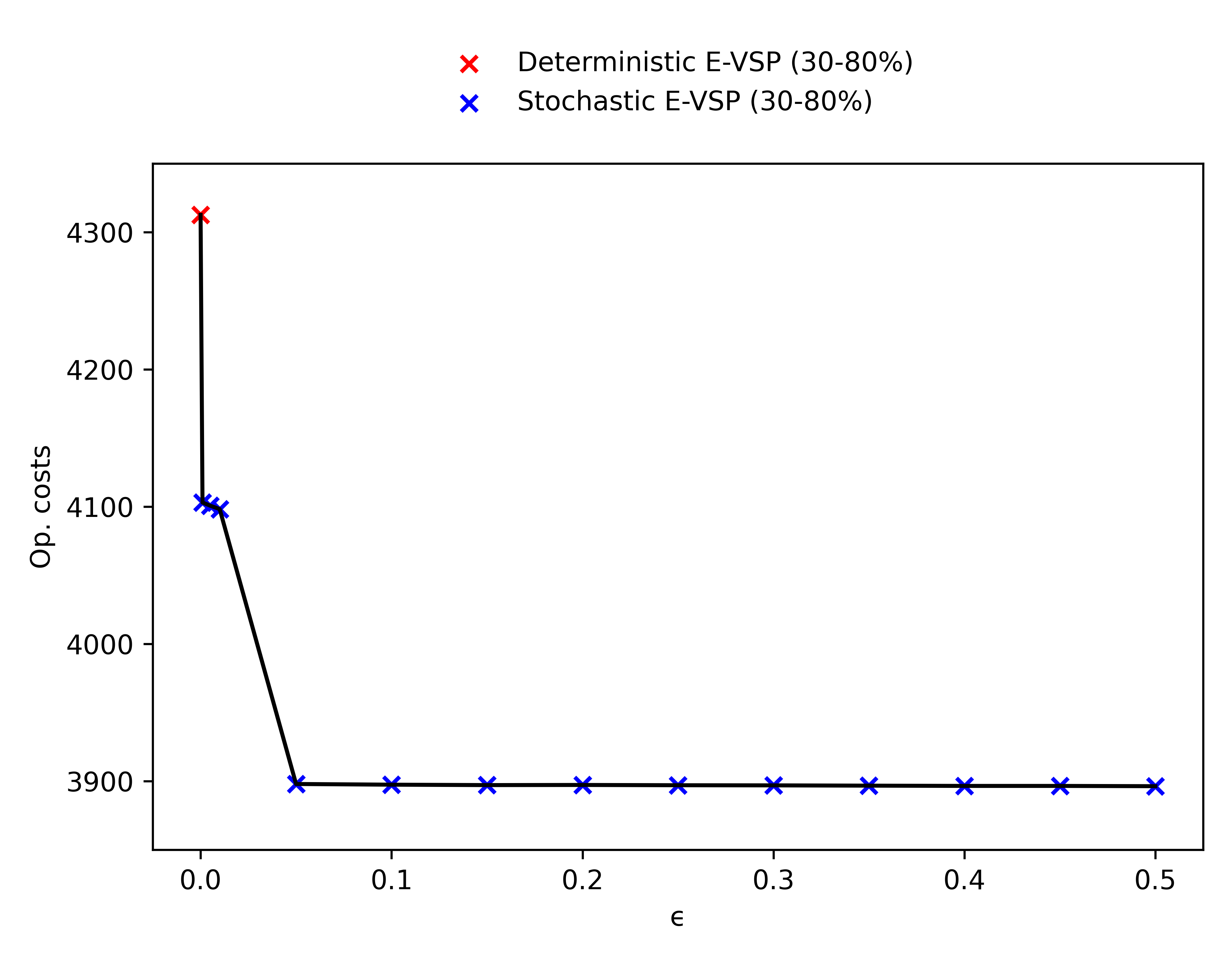}}
    \caption{\label{fig:results_I1}Comparison of the deterministic and the stochastic E-VSP with battery degradation solutions---instance family I1}
\end{figure}

\begin{figure}[ht!]
    \centering
    \subfloat[Range 20-80\%]{\label{fig:results_I2_20-80}\includegraphics[width=0.48\textwidth, height = 6.5 cm, keepaspectratio]{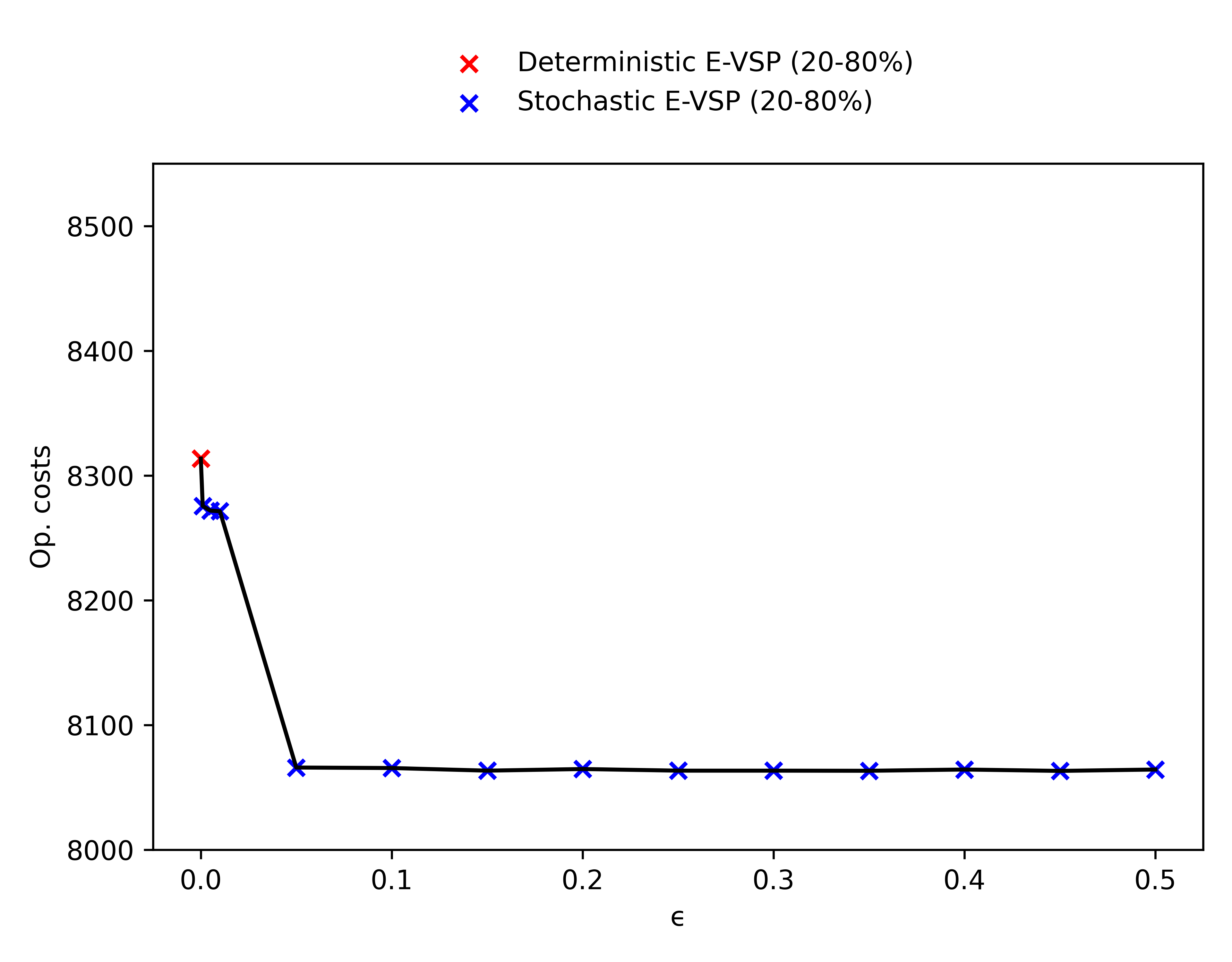}}\hspace{0cm}
    \subfloat[Range 30-80\%]{\label{fig:results_I2_30-30}\includegraphics[width=0.48\textwidth, height = 6.5 cm, keepaspectratio]{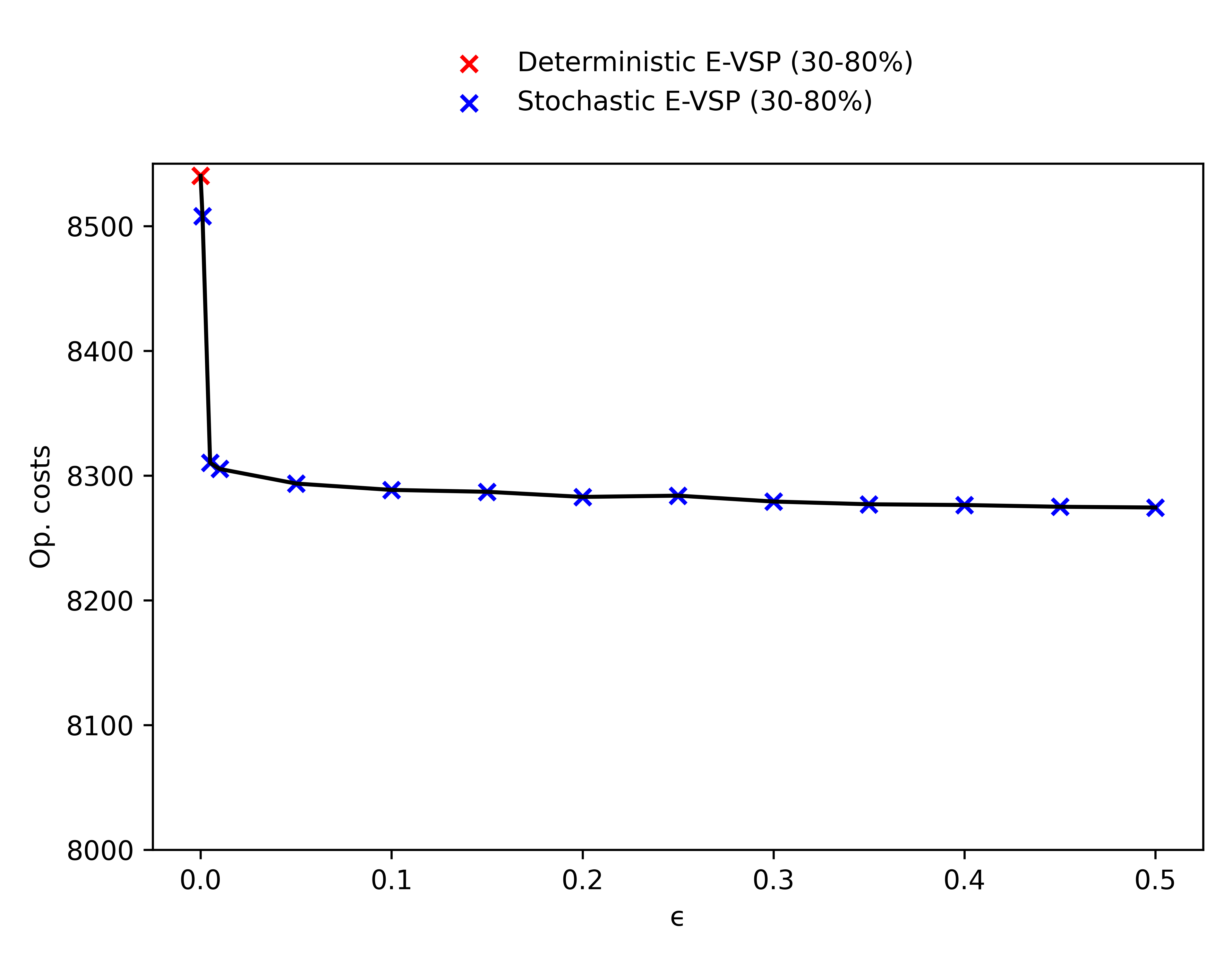}}
    \caption{\label{fig:results_I2}Comparison of the deterministic and the stochastic E-VSP with battery degradation solutions---instance family I2}
\end{figure}

\begin{figure}[ht!]
    \centering
    \subfloat[Range 20-80\%]{\label{fig:results_I3_20-80}\includegraphics[width=0.48\textwidth, height = 6.5 cm, keepaspectratio]{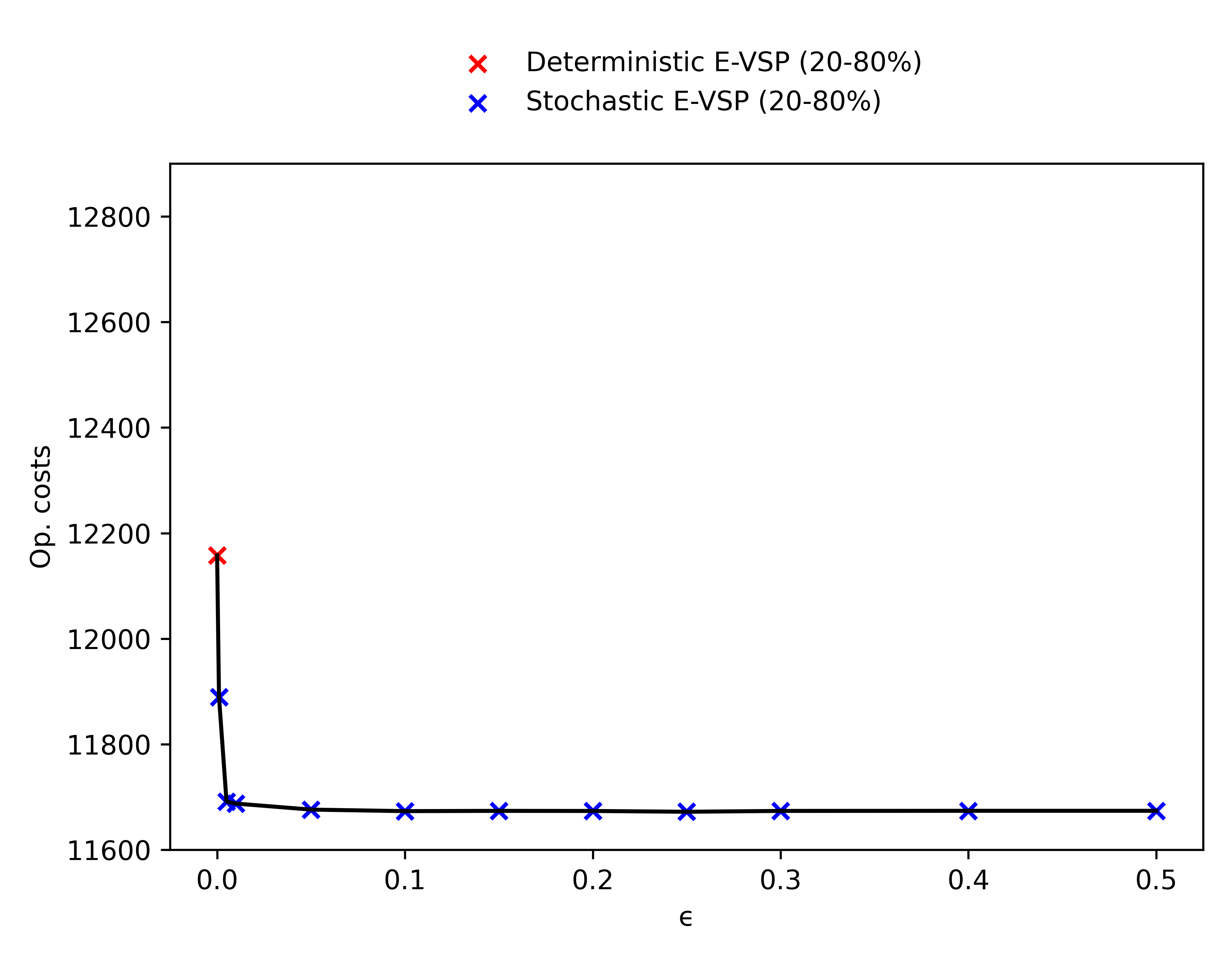}}\hspace{0cm}
    \subfloat[Range 30-80\%]{\label{fig:results_I3_30-80}\includegraphics[width=0.48\textwidth, height = 6.5 cm, keepaspectratio]{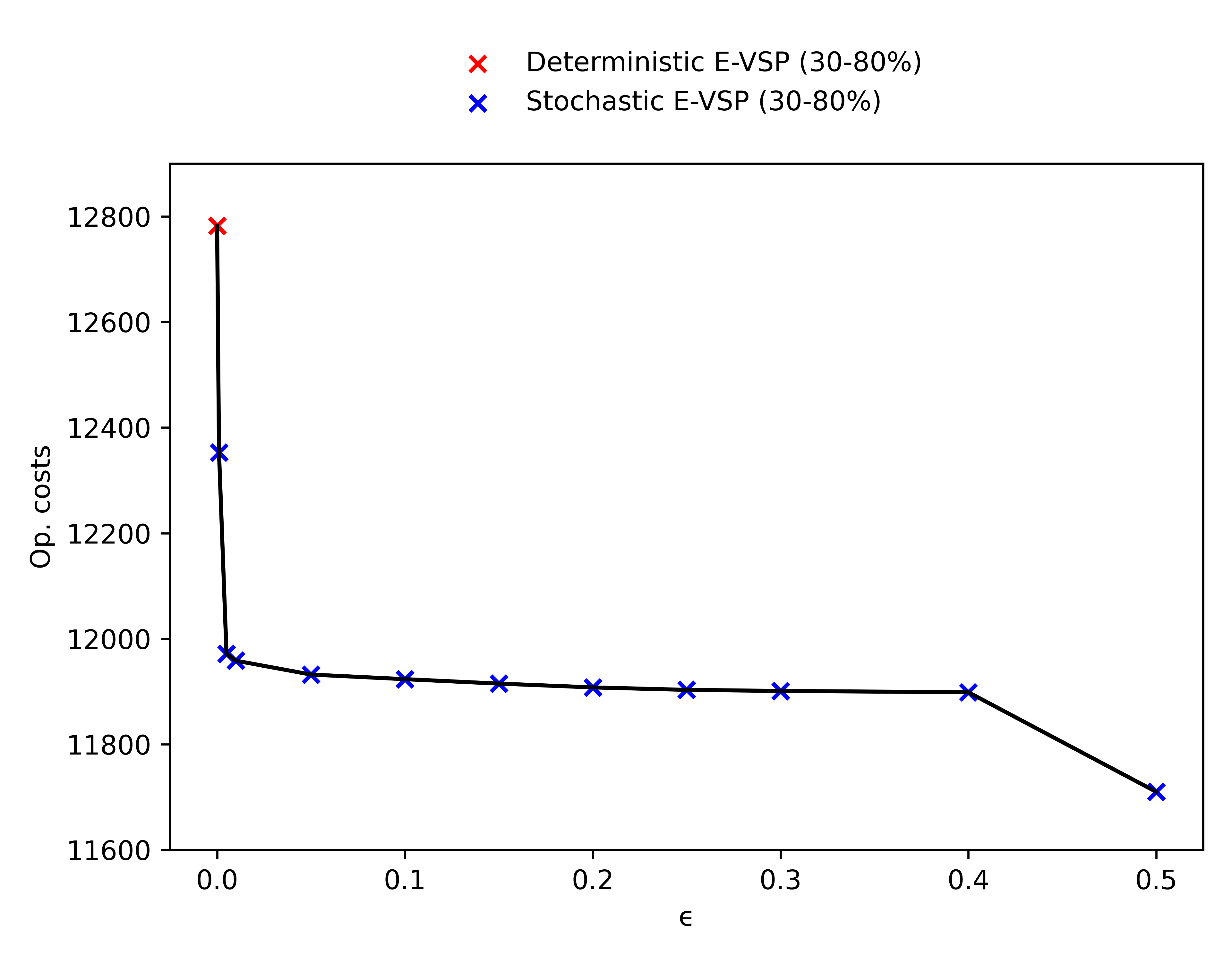}}
    \caption{\label{fig:results_I3}Comparison of the deterministic and the stochastic E-VSP with battery degradation solutions---instance family I3}
\end{figure}

Operational cost improvements generally stabilize after a certain threshold of \(\epsilon\), as shown in Figures \ref{fig:results_I1}-\ref{fig:results_I3}. For example, in I1, improvements plateau at 5.07\% for the 20–80\% range as \(\epsilon\) increases beyond 0.10. Similarly, for I2 and I3, cost improvements also stabilize, indicating diminishing returns from further relaxing the chance constraint. However, in I3, operational costs improve again for the 30–80\% range at \(\epsilon = 0.5\) (see Figure \ref{fig:results_I3_30-80}).

The operational costs are consistently higher for the 30–80\% SoC range compared to the 20–80\% range because the narrower operating range limits available energy resources. For instance, in I3 with \(\epsilon = 0.005\), the operational costs for the 20–80\% range are 11,691.6, compared to 11,971.1 for the 30–80\% range. However, the relative improvements in operational costs compared to the deterministic benchmark are generally larger for the 30–80\% range.

\subsection{Selecting the chance constraint's probability threshold}
\label{sec:selection_epsilon}

Selecting an appropriate probability threshold for the chance constraint is a critical decision in addressing the E-VSP with battery degradation and stochastic energy consumption. One effective approach involves employing a battery fading model to evaluate various solutions under different $\epsilon$ values. In this work, we use the battery fading model from \cite{Lam2012}. Note, however, that other capacity fading models can be employed in practice, or alternative methods, such as pilot programs, can be used to identify a suitable $\epsilon$ value.

The model from \cite{Lam2012} quantifies battery capacity fading based on the average SoC and deviation of the SoC during discharging and charging cycles. Given a vehicle schedule \(s \in \mathcal{S}^*\), we approximate the fading rate of the battery capacity by dividing \(s\) into discharging and charging cycles where energy consumption and replenishment occur. For simplicity, the SoC is assumed to decrease linearly with time. The capacity fading rate is first computed for each cycle and then aggregated for the entire schedule. We employ a Monte Carlo simulation with $K=1,000$ iterations to estimate the yearly capacity fading rate of a schedule, where at each iteration a random energy consumption value is generated for each trip in the schedule. Additional details on the methodology for computing battery degradation using the capacity fading model from \cite{Lam2012} can be found in Appendix \ref{section:appendixA}.

\begin{figure}
    \centering
    \includegraphics[width=0.7\linewidth]{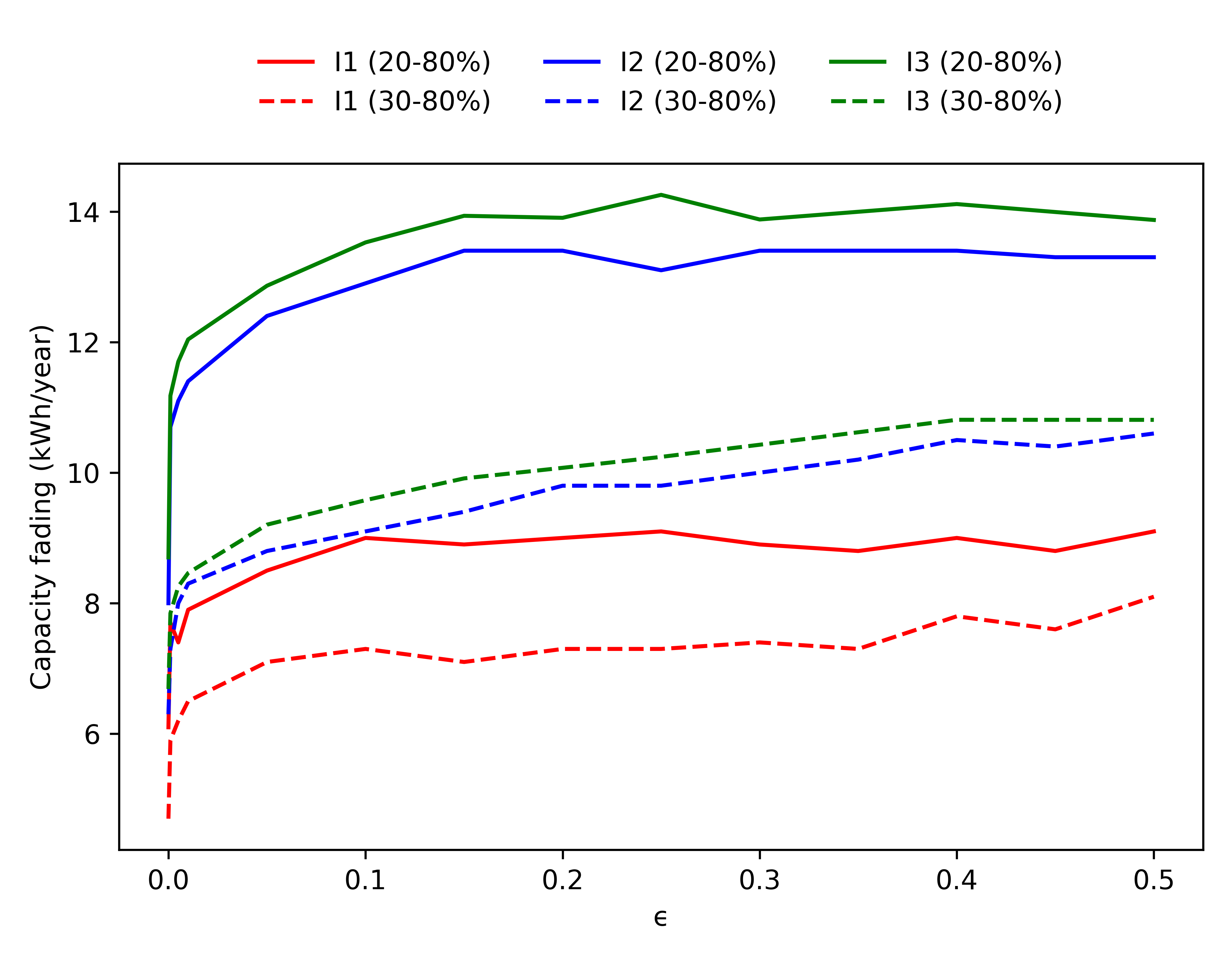}
    \caption{Average yearly capacity fading per vehicle vs. chance constraint probability threshold ($\epsilon$).}
    \label{fig:capacity_fading}
\end{figure}

Using this latter methodology, we analyze in Figure \ref{fig:capacity_fading} the relationship between the probability threshold ($\epsilon$) for overusing the batteries and the resulting average yearly capacity fading (measured in kWh/year) per vehicle. The results highlight a trade-off: higher $\epsilon$ values allow greater flexibility in operations but accelerate battery degradation. For instance, a capacity fading rate of 14 kWh/year---the highest observed for instance family I3 with a 20–80\% SoC range---corresponds to an expected battery lifespan of approximately 4.3 years. By comparison, a fading rate of 10 kWh/year translates to an estimated 6-year lifespan for a 300 kWh battery, assuming an end of life threshold of 80\%. However, it is important to note that these estimates may not fully reflect real-world conditions, as the model proposed by \cite{Lam2012} accounts only for cycle aging and neglects calendar aging. Despite these limitations, the presented curves provide insights for stakeholders to balance operational efficiency and battery longevity by selecting an appropriate $\epsilon$ value.

The analysis reveals varying degradation patterns across instance configurations (I1, I2, I3) and recommended SoC ranges (20–80\% and 30–80\%). Instance family I3 with a 20–80\% SoC range exhibits the highest capacity fading rate, while I1 with a 30–80\% range shows the lowest. Larger instance families (i.e., those with more timetabled trips) generally exhibit higher fading rates, as vehicles in these instances have to perform more trips per day on average---15.7 trips/EB for I1, 19.9 trips/EB for I2, and 21.8 trips/EB for I3 in the deterministic case with a 20–80\% SoC range. For all instance families, capacity fading rates increase steeply at low $\epsilon$ values but stabilize as $\epsilon$ exceeds 0.1. Nevertheless, the increase in fading rate with $\epsilon$ is less pronounced than the impact of changing the recommended SoC range. Across all instances, the 30–80\% SoC range consistently results in lower fading rates than the 20–80\% range.

Combining a 30–80\% SoC range with a small $\epsilon$ can yield win-win solutions compared to a deterministic approach. For instance, in instance family I3, an $\epsilon$ value of 0.005 with a 30–80\% SoC range reduces operational costs (11,971.1 vs. 12,157.7), the number of EBs required (11.2 vs. 11.4), and the capacity fading rate (8.25 vs. 8.70 kWh/year---a reduction of 5.2\%) compared to the deterministic benchmark with a 20–80\% range. This is also true for the instance family I2, where an $\epsilon$ value of 0.005 with a 30–80\% range reduces operational costs (8,301.6 vs. 8,311.9) and the capacity fading rate (7.96 vs. 8.05 kWh/year---a reduction of 1.1\%) compared to the deterministic benchmark with a 20–80\% range.

\section{Conclusions}
\label{sec:conclusions}

In this study, we presented a chance-constrained integer linear program for the E-VSP with battery degradation and stochastic energy consumption. The chance constraint ensures compliance with recommended SoC ranges, such as 20-80\% or 30-80\%, as specified by battery leasing companies or manufacturers, by limiting the probability of battery overuse. Unlike deterministic approaches that rely on worst-case assumptions, our method incorporates uncertainty in energy consumption, enabling more cost-efficient solutions. Additionally, the model avoids reliance on specific capacity fading functions, which are often calibrated using test data that fail to capture real-world driving conditions and are highly battery-specific. The model also accounts for partial en-route charging, nonlinear charging profiles, and limited-capacity charging stations. Nonlinear charging functions are approximated using piecewise linear functions.

To solve the model, we proposed a heuristic branch-and-price algorithm featuring stochastic pricing problems, specifically stochastic shortest path problems with resource constraints. A labeling algorithm with stochastic dominance rules based on the CDFs of SoC conditioned on the event that the battery remains within recommended SoC ranges is developed. Computing time is reduced by implementing a heuristic branching strategy, namely a diving heuristic, and a constraint perturbation method. 

Computational experiments conducted on realistic synthetic instances validated the performance of our heuristic, demonstrating its ability to derive near-optimal solutions. The results showed that small adjustments to the chance constraint probability threshold \(\epsilon\) can yield important cost savings. For instance, operational costs can be reduced by up to 10\% by allowing a 5\% probability of operating an EB outside the recommended SoC range. Using the battery fading function from \cite{Lam2012}, we further investigated the impact of \(\epsilon\) on battery degradation and showed that our approach can provide win-win solutions, reducing both operational costs and battery wear compared to deterministic baselines.

While our work demonstrates the added value of our stochastic model for the E-VSP with battery degradation, it remains limited to relatively small instances. Future research could focus on reducing computational times to enable the solution of larger instances. 

\section{Acknowledgments}
We thank the Optimization Algorithm Development team at GIRO Inc. for the insightful discussions that helped shape the problem definition. This work was partially supported by GIRO Inc. and the Natural Sciences and Engineering Research Council of Canada (NSERC) under grant RDCPJ 520349-17. Léa Ricard gratefully acknowledges the financial support of NSERC through the BESC D3-558645-2021 scholarship. The work of the first and third authors was supported by the Canada Excellence Research Chair in  `Data Science for Real-Time Decision-Making'.

\bibliographystyle{bib_style}
\bibliography{refs}

\newpage
\appendix

\section{Computing the battery degradation level of an EB schedule}
\label{section:appendixA}
We evaluate the degradation level of the battery of an EB associated with a schedule $s$ by analyzing the capacity fading encountered during each discharging and charging cycle $\psi_t \in s$. A cycle involves two primary stages: 1) discharging, which occurs when the vehicle is in use and draws power from the battery to propel itself, and 2) charging, which happens when the battery is connected to an external power source to restore its energy. Specifically, let us consider a vehicle schedule $s = (n_o^d, v_1, h^{w1}_{r0}, h^{c1}_{r1}, h^{w1}_{r1}, v_2, v_3, h^{w1}_{r3},h^{c1}_{r4},h^{w1}_{r4},v_6,n_1^d)$ in the network illustrated in Figure \ref{fig:network}. This vehicle schedule consists of three discharging and charging cycles: $\psi_1 = (n_o^d, v_1, h^{w1}_{r0}, h^{c1}_{r1}, h^{w1}_{r1})$, $\psi_2 = (v_2, v_3, h^{w1}_{r3},h^{c1}_{r4},h^{w1}_{r4})$, and $\psi_3 = (v_6,n_1^d)$.

In this paper, we use the capacity fading function developed in \cite{Lam2012} to approximate the degradation level induced by a solution to the E-VSP with battery degradation and stochastic energy consumption. This model quantifies the rate at which a battery's capacity degrades, i.e., the loss in kWh of the battery's capacity per 1 kWh of charge processed. The degradation rate depends on both the average SoC and the SoC deviation.  However, in real-world scenarios, the discharging and charging cycles differ daily due to the randomness of energy consumption. Consequently, the average SoC and SoC deviation also fluctuate. To address this, we employ a stochastic simulation approach. In Section \ref{sec:capacity_fading_function}, we present our adaptation of the function of \cite{Lam2012} to meet our specific scenario. Additionally, we introduce in Section \ref{sec:Simul} a Monte Carlo simulation that approximates the expected capacity fading of the battery of an EB over the entire planning horizon (typically one day).

\subsection{Battery capacity fading function}
\label{sec:capacity_fading_function}

Let $\sigma_{\psi_t}^b$, $\sigma_{\psi_t}^m$, and $\sigma_{\psi_t}^e$ be the SoC at the beginning, middle (after the discharging regime), and end (after recharging) of the discharging and charging cycle $\psi_t$, respectively. The average SoC in cycle $\psi_t$ is defined as

\begin{equation}
\label{eq:avg_SoC}
    \sigma_{\psi_t}^{avg} = \frac{\sigma_{\psi_t}^b + \sigma_{\psi_t}^m + \sigma_{\psi_t}^e}{3},
\end{equation}

\noindent and the SoC deviation is defined as 

\begin{equation}
\label{eq:SoC_dev}
    \sigma_{\psi_t}^{dev} = \sigma_{\psi_t}^{avg} - \sigma_{\psi_t}^m.
\end{equation}

Note that equations (\ref{eq:avg_SoC}) and (\ref{eq:SoC_dev}) assume for simplicity that the charging and discharging processes are linear, even though the SoC varies non-linearly with time. The capacity fading rate (in kWh/1 kWh processed) can be computed as 

\begin{equation}
    \phi(\sigma^b_{\psi_t}, \sigma^m_{\psi_t}, \sigma^e_{\psi_t}) = \gamma_1 \sigma_{\psi_t}^{dev} e^{\gamma_2 \sigma_{\psi_t}^{avg}} + \gamma_3 e^{\gamma_4 \sigma_{\psi_t}^{dev}},
\end{equation}

\noindent where $\gamma_1, \dots, \gamma_4$ are model parameters borrowed from \cite{Lam2012} which are set as: $\gamma_1 =$ -4.092e-4, $\gamma_2 = -2.167, \gamma_3 =$ 1.408e-5, and $\gamma_4 = 6.130$.

Given a charge processed of $\Sigma_{\psi_t} =(\sigma_{\psi_t}^b -  \sigma_{\psi_t}^m) + (\sigma_{\psi_t}^e - \sigma_{\psi_t}^m)$ in cycle $\psi_t$, the daily capacity fade of cycle $\psi_t$ (in kWh) can be computed as

\begin{equation}
     \Phi(\sigma^b_{\psi_t}, \sigma^m_{\psi_t}, \sigma^e_{\psi_t}) = \Sigma_{\psi_t} \times \phi (\sigma^b_{\psi_t}, \sigma^m_{\psi_t}, \sigma^e_{\psi_t}).
\end{equation}

The daily capacity fade of schedule $s$ is obtained by summing the daily capacity fade of every cycle $\psi_t \in s$ as 

\begin{equation}
    Q_{s} =  \sum_{\psi_t \in s} \Phi(\sigma^b_{\psi_t}, \sigma^m_{\psi_t}, \sigma^e_{\psi_t}),
\end{equation}

\noindent where we assume the SoC at the end of the last cycle is $\sigma^{init}$, i.e., the EBs are recharged overnight at the depot.

\subsection{Monte Carlo simulation}
\label{sec:Simul}

This section describes a Monte Carlo simulation to approximate offline the daily capacity fading yielded by a solution to the E-VSP with battery degradation and stochastic energy consumption. The Monte Carlo simulation is presented in Algorithm \ref{al:monte_carlo} and can be summarized as follows. For each vehicle schedule $s \in \mathcal{S}^*$ with $q_s$ cycles and in each iteration $k = 1, \dots, K$, a capacity fading value $Q_s^k$ is computed by iterating through all the nodes in the discharging and charging cycles $\psi_t \in  s \setminus \{\psi_{q_s}\}$ in Steps 7-24 and the last cycle $\psi_{q_s}$ in Steps 26-29. The capacity fading $\phi(\sigma^{b,k}_{\psi_t}, \sigma^{m,k}_{\psi_t}, \sigma^{e,k}_{\psi_t})$ of each cycle $\psi_t \in  s \setminus \{\psi_{q_s}\}$ is computed and added to $Q_s^k$ in Step 22. The same is done for the last cycle in Step 31. To do this, random energy consumption values $\mu_i^k$ are generated for each timetabled trip $i \in s$ in Steps 11 and 27. After iterating through all cycles in $s$, the algorithm adds the capacity fading of schedule $s$ in iteration $k$, $Q_s^k$, to the total capacity fading ($W$) in Step 32. This process is repeated $K$ times for each vehicle schedule $s \in \mathcal{S}^*$. Finally, in Step 35, the algorithm returns the approximated daily capacity fading per vehicle of the E-VSP solution.

\begin{algorithm}[!ht]
\newcommand\mycommfont[1]{\ttfamily\textcolor{blue}{#1}}
\SetCommentSty{mycommfont}
\DontPrintSemicolon

\SetAlgoLined
\SetKwInOut{Input}{input}
\SetKwInOut{Output}{output}
\SetKwFunction{AND}{and}

$W \leftarrow 0$\;
\For{$s \in \mathcal{S}^{*}$}{
$Q_s \leftarrow 0$\;
\For{$k\leftarrow 1$ \KwTo $K$}{
$\sigma^{b,k}_{\psi_1} \leftarrow \sigma^{\text{init}}$\;
$Q_s^k \leftarrow 0$\;
\For{$\psi_t \in s \setminus \{\psi_{q_s}\}$}{
$\sigma^{m,k}_{\psi_t} \leftarrow \sigma^{b,k}_{\psi_t}$\;
\For{$i \in \psi_t$}{
\If{$i \in \mathcal{V}$}{
Randomly generate  $\mu_i^k$ from $e_i(\mu)$\;
$\sigma^{m,k}_{\psi_t} \leftarrow \sigma^{m,k}_{\psi_t}  - \mu_i^k - \iota_{i-1,i}$\;
}
\If{$i \in \mathcal{H}^W \text{ and } i-1 \in \mathcal{V}$}{
$\sigma^{m,k}_{\psi_t} \leftarrow \sigma^{m,k}_{\psi_t} - \iota_{i-1,i}$\;
$\sigma^{e,k}_{\psi_t} \leftarrow \sigma^{m,k}_{\psi_t}$\;
}
\If{$i \in \mathcal{H}^C$}{
$\sigma^{e,k}_{\psi_t} \leftarrow \lambda(\sigma^{e,k}_{\psi_t}, \rho)$\;
}
}
$Q_s \leftarrow Q_s^k + \phi(\sigma^{b,k}_{\psi_t}, \sigma^{m,k}_{\psi_t}, \sigma^{e,k}_{\psi_t})$\;
$\sigma^{b,k}_{\psi_{t+1}} \leftarrow \sigma^{e,k}_{\psi_t}$\;
}
$\sigma^{m,k}_{\psi_{q}} \leftarrow \sigma^{e,k}_{\psi_{t}}$\;
\For{$i \in \psi_{q_s}$}{
Randomly generate  $\mu_i^k$ from $e_i(\mu)$\;
$\sigma^{m,k}_{\psi_{q_s}} \leftarrow \sigma^{m,k}_{\psi_{q_s}} - \mu_i^k - \iota_{i-1,i}$\;
}
$\sigma^{e,k}_{\psi_{q_s }} \leftarrow \sigma^{init}$\;
$Q_s \leftarrow Q_s^k + \phi(\sigma^{b,k}_{\psi_{q_s}}, \sigma^{m,k}_{\psi_{q_s}}, \sigma^{e,k}_{\psi_{q_s }})$\;
$W \leftarrow W + Q_s^k $\;
}
}
Return $W/(K \times |\mathcal{S}^*|)$\;
\caption{Simulation to approximate the total daily capacity fade induced by a E-VSP solution}
\label{al:monte_carlo}
\end{algorithm}

\end{document}